\documentclass[a4paper,leqno,12pt]{amsart}
\usepackage{amssymb,amscd}
\usepackage{graphicx}
\usepackage{xy}
\usepackage{fullpage}
\xyoption{all}
\divide\oddsidemargin by 2
\usepackage{enumerate}
\makeatletter
\renewcommand\theenumi{\@alph\c@enumi}
\renewcommand\theenumii{\@alph\c@enumii}
\renewcommand\theenumiii{\@alph\c@enumiii}
\renewcommand\theenumiv{\@alph\c@enumiv}

\let\@listlla\list

\def\list#1#2{\@listlla{#1}{#2\itemsep=2pt\parsep=0pt\topsep=3pt plus 1pt minus 1 pt}}
\newcommand{\tridiagram}[6]{{\par\par \centering
\@picture(120,120)(0,0)
\put(30,95){\makebox(0,0)[r]{$#1$}}
\put(90,30){\makebox(0,0)[tl]{$#3$}}
\put(90,95){\makebox(0,0)[l]{$#2$}}
\put(60,102){\makebox(0,0)[b]{$#4$}}
\put(102,60){\makebox(0,0)[l]{$#6$}}
\put(50,50){\makebox(0,0)[tr]{$#5$}}
\thinlines
\put(40,95){\vector(1,0){40}}
\put(95,80){\vector(0,-1){40}}
\put(25,80){\vector(1,-1){55}}
\endpicture\par\par}\noindent\ignorespaces}
\def\@map#1#2[#3]{\mbox{$#1 \colon #2 \longrightarrow #3$}}
\def\map#1#2{\@ifnextchar [{\@map{#1}{#2}}{\@map{#1}{#2}[#2]}}

\newcommand{\vois}[3]{{\mathcal N}^{#1}_{#2}(#3)}
\newcommand{\mz}{\mathbb{Z}}

\newcommand{\F}[1]{\ensuremath{\mathbb{F}_{#1}}}

\makeatother
\newtheorem{theorem}{Theorem}[section]
\newtheorem{proposition}[theorem]{Proposition}
\newtheorem{corollary}[theorem]{Corollary}
\newtheorem{lemma}[theorem]{Lemma}
\theoremstyle{definition}
\newtheorem{definition}[theorem]{Definition}

\newtheorem{remark}[theorem]{Remark}
\newcommand{\pol}[2]{{\noindent \small \bf Proof of Lemma
\ref{#1}: }{\rm #2 \hfill $\Box$}}
\newcommand{\pop}[2]{{\noindent \small \bf Proof of Proposition
\ref{#1}: }{\rm #2 \hfill $\Box$}}

\newcommand{\pot}[2]{{\noindent \small \bf Proof of Theorem
\ref{#1}: }{\rm #2 \hfill $\Box$}}

\newcommand{\pols}[1]{{\noindent \small \bf Proof: }{\rm #1 \hfill $\Box$}}

\newcommand{\pocs}[1]{{\noindent \small \bf Proof: }{\rm #1 \hfill $\Box$}}
\newcommand{\pots}[1]{{\noindent \small \bf Proof: }{\rm #1 \hfill $\Box$}}
\begin{document}

\title{Geodesics in trees of hyperbolic and relatively hyperbolic
groups}

\author[F. Gautero]{Fran\c{c}ois Gautero}
\address{Universit\'e Blaise Pascal,
Campus des C\'ezeaux, Laboratoire de Math\'ematiques, 63177
Aubi\`ere, France}

\email{Francois.Gautero@math.univ-bpclermont.fr}

\thanks{{\em Acknowledgements.} This paper benefited
from discussions of the author with F. Dahmani, V. Guirardel, M.
Heusener, I. Kapovich and M. Lustig. Warm thanks are particularly
due to I. Kapovich for his explanations about Bestvina-Feighn paper.
At that time, the author was Assistant at the University of Geneva,
whereas I. Kapovich visited this university thanks to a funding of
the Swiss National Science Foundation. The author is glad to
acknowledge support of these institutions. Finally, P. de la Harpe
also deserves a great share of these acknowledgements for his help.}



%
%

\begin{abstract}
We present a careful approximation of the geodesics in trees of
hyperbolic or relatively hyperbolic groups. As an application we
prove a combination theorem for finite graphs of relatively
hyperbolic groups, with both Farb's and Gromov's definitions.
\end{abstract}

\maketitle

%
%

\section{Introduction}

The main part of this paper is devoted to give a precise description
of the geodesics in trees of hyperbolic (and thereafter relatively
hyperbolic - see below) groups. Such a work might appear not very
appealing, and somehow quite technic. In order to show that this
however might be worthy, let us give an application: a {\em
combination theorem} for hyperbolic and relatively hyperbolic
groups. That is, a theorem giving a condition for the fundamental
group of a graph of relatively hyperbolic groups being a relatively
hyperbolic group. In \cite{BF} (see also \cite{IKa}), the authors
introduce the notion of (finite) qi-embedded graph of groups and
spaces $\mathcal G$. Then, assuming the Gromov hyperbolicity of the
vertex spaces and the quasiconvexity of the edge spaces in the
vertex spaces, they give a criterion for the hyperbolicity of the
fundamental group of $\mathcal G$. Since then different proofs have
appeared, which treat the so-called `acylindrical case': see, among
others, \cite{Gi,KM}. Acylindrical means that the fixed set of the
action of any element of the fundamental group of the graph of
groups on the universal covering has uniformly bounded diameter. The
non-acylindrical case is less common: see \cite{Ika1} which relies
on \cite{BF} but clarifies its consequences when dealing with a
certain class of mapping-tori of injective, non surjective free
group endomorphisms, or \cite{Ga} which, by an approach similar to
the one presented here, gives a new proof of \cite{BF} in the case
of mapping-tori of free group endomorphisms. Nowadays the attention
has drifted from hyperbolic groups to {\em relatively hyperbolic
groups}. A notion of relative hyperbolicity was already defined by
Gromov in his seminal paper \cite{Gr}. Since then it has been
revisited and elaborated on in many papers. Two distinct definitions
now coexist. In parallel to the Gromov relative hyperbolicity,
sometimes called {\em strong relative hyperbolicity}, there is the
sometimes called {\em weak relative hyperbolicity} introduced by
Farb \cite{Fa}. Bowditch \cite{Bo} and Osin \cite{Os} give
alternative definitions, but which are equivalent either to Farb's
or to Gromov's definition. In fact, it has been proved \cite{Bu,Os}
(also \cite{Bo}) that Gromov definition is equivalent to Farb
definition plus an additional property termed Bounded Coset
Penetration property (BCP in short), due to Farb \cite{Fa}.
Relatively hyperbolic groups in the strong (that is Gromov) sense
form a class encompassing hyperbolic groups, fundamental groups of
geometrically finite orbifolds with pinched negative curvature,
groups acting on CAT(0)-spaces with isolated flats among many
others. First combination theorems in some particular (essentially
acylindrical) cases have been given in the setting of the relative
hyperbolicity: \cite{Ali}, \cite{Da} or \cite{Os1,Os2}. One result
\cite{Ga2} treats a particular non-acylindrical case, namely the
relative hyperbolicity of one-ended hyperbolic by cyclic groups.
Since a first version of this paper was written, a paper \cite{MR}
appeared on the arXiv, also giving a combination theorem dealing
with more general non-acylindrical cases than \cite{Ga2}: the
authors heavily rely upon \cite{BF}, which they use as a
``black-box''. Getting a ``general'' combination theorem for
relatively hyperbolic groups is one of the questions (attributed to
Swarup) raised in Bestvina's list \cite{Bliste}. We offer here an
answer, as an application of our work on geodesics in trees of
spaces. We would like to emphasize at once that we do not appeal to
\cite{BF}, but instead give a new proof of it as a particular case.
Where the authors of \cite{BF} use ``second-order'' geometric
characterization of hyperbolicity via isoperimetric inequalities, we
use ``first-order'' geometric characterization, via approximations
of geodesics and the thin triangle property. At the expense of
heavier and sometimes tedious computations, this na\"{\i}ve approach
allows us to engulf in a same setting (at least when dealing with
combination theorems) both absolute and relative hyperbolicity.

In order to illustrate our results, let us just give the following
particular case:

\begin{theorem}
\label{theoreme principal} Let $G$ be a finitely generated group
and let $\mathcal H$ be
 a finite family of subgroups  of $G$. Let $\F{r}$ be a
uniform free group of relatively hyperbolic automorphisms of $(G,
{\mathcal H})$. If $G$ is weakly hyperbolic relative to $\mathcal
H$, then $G \rtimes \F{r}$
 is weakly hyperbolic relative to $\mathcal
 H$. If $G$
is strongly hyperbolic relative to $\mathcal H$, then $G \rtimes
\F{r}$
 is strongly hyperbolic relative to a $\F{r}$-extension of $\mathcal H$.
\end{theorem}

See \ref{premieres definitions}, \ref{auto relativement
hyperbolique} and \ref{extension} for the definitions of (free
groups of) relatively hyperbolic automorphisms and of a free
extension of a family. Theorem \ref{theoreme principal} is a
compilation of Theorems \ref{theoreme principal faible} and
\ref{theoreme principal fort}. When $r=1$ in the above theorem, that
is when the considered free group is just $\mz$, we get the
classical ``mapping-torus'' case, that is the case of semi-direct
products $G \rtimes \mz$ with $G$ a relatively hyperbolic group.
Corollary \ref{rappel0} gives a concrete application, when $G$ is
the fundamental group of a compact surface and $\mz$ acts on $G$ by
an automorphism induced by a homeomorphism $h$ of the surface $S$.
In this case the mapping-torus group $G \rtimes \mz$ is weakly
hyperbolic relative to the family formed by the cyclic subgroups
generated by the boundary loops of $S$, the subgroups defined (up to
conjugacy) by the maximal subsurfaces of $S$ preserved up to isotopy
by (a power of) the homeomorphism, and the cyclic subgroups
generated by the reduction curves which are not already contained in
the previous subgroups. It is strongly hyperbolic relatively to the
family of subgroups composed of the subgroups associated to the
boundary tori (assume for simplicity that $S$ is orientable and that
$h$ preserves the orientation), the subgroups associated to the
$3$-dimensional submanifolds which are the mapping-tori of the
maximal non pseudo-Anosov components and the subgroups associated to
the $2$-dimensional tori which are the mapping-tori of the remaining
reduction curves.

Theorem \ref{theoreme principal} is only a particular, but
emblematic, case of our more general results, see Theorems
\ref{intro sup}, \ref{theoreme principal faible} and Corollary
\ref{premier corollaire faible} in Section \ref{rh} and Theorems
\ref{je suis content de moi 2}, \ref{theoreme principal fort},
\ref{injectif} and Corollaries \ref{premier corollaire fort} and
\ref{rappel0} in Section \ref{strong}. All are consequences of
Theorems \ref{pipeau} and \ref{pipeau50} in Sections \ref{aqgpc} and
\ref{aqggc} about the behavior of quasi geodesics in trees of
hyperbolic spaces.

\subsection{Plan of the paper:} Section \ref{basis} contains the basis,
from quasi isometries to the ``hallways-flare'' property. Section
\ref{aqgpc} deals with the approximation of quasi geodesics in the
particular case where all the attaching-maps of the considered tree
of hyperbolic spaces are quasi isometries. Section \ref{aqggc}
contains the adaptations to the general case. The important notions
appearing in these two sections are the corridors in Section
\ref{aqgpc}, and the generalized and pseudo corridors in Section
\ref{aqggc}. These two sections appeal to three important
Propositions whose proofs are delayed: Proposition \ref{proposition
importante 1} is proved in Section \ref{propimpo1}, Proposition
\ref{quasiconvexite} in Section \ref{quasi convexe} and Proposition
\ref{qc} in subsection \ref{popqc} of this last section. Section
\ref{rh} presents the results about the hyperbolicity and the weak
relative hyperbolicity whereas Section \ref{strong} deals with the
consequences about the strong relative hyperbolicity of graphs of
strongly relatively hyperbolic groups. This last section contains
another proposition whose proof is postponed for a while, to
subsection \ref{upc} of Section \ref{quasi convexe}.

\section{Preliminaries}

\label{basis}

If  $(X,d)$ is a metric space with distance function $d$, and $x$ a
point in $X$, we set $B_{x}(r) = \{y \in X \mbox{ ; } d(x, y) \leq
r\}$. If $A$ and $B$ are any two subsets of $(X,d)$, $d^{i}(A,B) =
\displaystyle \inf_{x \in A,y \in B} d(x,y)$. We set also
$\vois{r}{d}{A} = \{x \in X \mbox{ ; } d^{i}(x,A) \leq r\}$ and
 $d^{H}(A, B) = \sup \{r \geq 0 \mbox{ ; } A \subset
 \vois{r}{d}{B} \mbox{ and } B \subset \vois{r}{d}{A}\}$ is then the
 usual Hausdorff distance between  $A$ and $B$. Finally,
 $\mathrm{diam}_X(A)$ stands for $\displaystyle \sup \{d(x,y) \mbox{ ; }  (x,y) \in A \times
 A \}$.

\subsection{Quasi isometries, quasi geodesics and hyperbolic spaces}

A {\em $(\lambda,\mu)$-quasi isometric embedding} from $(X_{1},
d_{1})$ to $(X_{2}, d_{2})$ is a map  $f \colon X_1 \rightarrow
X_2$ such that, for any $x,y$ in $X_{1}$:

 $$ \frac{1}{\lambda} d_{1}(x, y) - \mu \leq d_{2}(f(x), f(y)) \leq
 \lambda d_{1}(x, y) + \mu $$

 A {\em $(\lambda,\mu)$-quasi isometry}
 $f \colon (X_1,d_1) \rightarrow (X_2,d_2)$
 is a $(\lambda,\mu)$-quasi isometric embedding such that for any
 $y \in X_2$ there exists $x \in X_1$ with $d_2(f(x),y) \leq \mu$.

 A
{\em $(\lambda,\mu)$-quasi geodesic} in a metric space $(X,d)$ is
the image of an interval of the real line under a
$(\lambda,\mu)$-quasi isometric embedding.

We work with a version of the Gromov hyperbolic spaces which is
slightly extended with respect to the most commonly used. We do not
require first that they be geodesic, and second that they be proper,
that is closed balls are not necessarily compact. Instead of
geodesicity, we require {\em quasi-geodesicity}. We say that a
metric space $(X,d)$ is a {\em $(r,s)$-quasi geodesic space} if, for
any two points $x,y$ in $X$ there is a $(r,s)$-quasi geodesic
between $x$ and $y$. We then denote by $[x,y]$ such a $(r,s)$-quasi
geodesic (and of course in a geodesic space, $[x,y]$ denotes any
geodesic between $x$ and $y$). A quasi geodesic metric space is a
metric space which is $(r,s)$-quasi geodesic for some non negative
real constants $r,s$. The $(r,s)$-quasi geodesic triangles in a
$(r,s)$-quasi geodesic metric space $(X, d)$ are {\em thin} if there
exists $\delta \geq 0$ such that any $(r,s)$-quasi geodesic triangle
in $(X, d)$ is {\em $\delta$-thin}, that is any side is contained in
the $\delta$-neighborhood of the union of the two other sides. In
this case, $X$ is a {\em $\delta$-hyperbolic space}. A metric space
$(X,d)$ is a {\em Gromov hyperbolic space} if there exists $\delta
\geq 0$ such that $(X,d)$ is a $\delta$-hyperbolic space. The slight
``generalization'' from geodesic to quasi geodesic spaces is only a
technical point. But not requiring our spaces to be proper is
important in order to deal with relatively hyperbolic groups, the
definitions of which involve non-proper metric graphs.

\subsection{Trees of spaces}

A {\em metric tree} is a simplicial tree with all edges isometric to
$(0,1)$. If $\mathcal T$ is a metric tree, we denote by
$|.|_{{\mathcal T}}$ the length of a path in $\mathcal T$ and by
$d_{\mathcal T}$ the associated distance.

\begin{definition} (compare \cite{BF})
\label{coqi} A {\em tree of metric spaces $(\tilde{X},{\mathcal
T},\pi)$} is a metric space $\tilde{X}$ equipped with a projection
$\pi \colon \tilde{X} \rightarrow \mathcal T$ onto a metric tree
$\mathcal T$ which satisfy the following properties for some
$\lambda \geq 1$ and $\mu \geq 0$:

\begin{enumerate}
  \item If $m_e$ is the midpoint of the edge $e$, then
  $\pi^{-1}(m_e) = X_e$ is a geodesic metric space and $\pi^{-1}(e)$ is
  isometric to $X_e \times (0,1)$.
  \item If $v$ is a vertex of $\mathcal T$, if ${\mathcal T}_s$
  is the tree $\mathcal T$ subdivided at the midpoints of the
  edges and $S_v$ is the closed
  star of $v$ in ${\mathcal T}_s$, then:
  \begin{itemize}
    \item $\pi^{-1}(v)$ is a geodesic metric space $X_v$;
    \item $\pi^{-1}(S_v)$ is obtained from the disjoint union of $X_v$ with the
    spaces $X_e \times [0,1/2]$, $e$ the edges of ${\mathcal T}_s$ in
    $S_v$, by identifying each $X_e \times \{0\}$ with a subset of
    $X_v$ under a $(\lambda,\mu)$-quasi isometric embedding.
  \end{itemize}
\end{enumerate}

A set $\pi^{-1}(x)$, $x \in {\mathcal T}$, is a {\em stratum}.

A {\em tree of hyperbolic spaces} is a tree of metric spaces such
that there is $\delta \geq 0$ for which the strata are
$\delta$-hyperbolic spaces.
\end{definition}

By definition, each stratum in a tree of metric spaces comes with a
distance, termed {\em horizontal distance}. A path contained in a
stratum is a {\em horizontal path} and we will also speak of the
{\em horizontal length} of a horizontal path.

\subsection{The telescopic metric}

\begin{definition}
Let $(\tilde{X},{\mathcal T},\pi)$ be a tree of metric spaces.

A {\em $v$-vertical segment} (resp. {\em $v$-vertical tree}) in
$\tilde{X}$ is (the image of) a section $\sigma_w$ (resp.
$\sigma_T$) of $\pi$ over a geodesic $w$ of $\mathcal T$ (resp. over
a subtree $T$ of $\mathcal T$) which is a $(v,v)$-quasi isometric
embedding.

The $\mathcal T$-length $|w|_{\mathcal T}$ is the {\em vertical
length} of the $v$-vertical segment $\sigma_w \colon w \rightarrow
\tilde{X}$.

If $x$ is a point in $\tilde{X}$ and $w$ is a geodesic of $\mathcal
T$ starting at $\pi(x)$, the notation $w x$ will denote the set of
points $y \in \tilde{X}$ such that some $v$-vertical segment $s$
with $\pi(s) = w$ connects $x$ to $y$ (in particular any such $y$
belongs to $\pi^{-1}(t(w))$).
\end{definition}

\begin{definition}
\label{ciel} Let $(\tilde{X},{\mathcal T},\pi)$ be a tree of metric
spaces.

A {\em $v$-telescopic path} is a path $p$ in $\tilde{X}$ which
satisfies the following properties:

\begin{itemize}
  \item $\pi(p)$ is an edge-path between
  two vertices of $\mathcal T$,
  \item $p$ is the
    concatenation of
horizontal paths in the strata over the vertices of $\mathcal T$
  and of non-trivial $v$-vertical segments.
  \end{itemize}

\end{definition}

\begin{definition}
\label{television} Let $p$ be a $v$-telescopic path in a tree of
metric spaces $(\tilde{X},{\mathcal T},\pi)$.

\begin{enumerate}
  \item  The {\em vertical
length} $|p|^v_{vert}$ of $p$ is the sum of the vertical lengths of
the maximal $v$-vertical segments. The {\em horizontal length}
$|p|^v_{hor}$ is the sum of the horizontal lengths of the maximal
horizontal subpaths in the complement of the maximal $v$-vertical
segments.
  \item The {\em telescopic length $|p|^{v}_{tel}$} of a
$v$-telescopic path $p$ is the sum of its horizontal and vertical
lengths.
\end{enumerate}
\end{definition}

\begin{definition}
\label{definition distance} Let $(\tilde{X},{\mathcal T},\pi)$ be a
tree of metric spaces. The {\em $v$-telescopic distance
$d^{v}_{tel}(x,y)$} between two points $x$ and $y$ is the infimum of
the telescopic lengths of the $v$-telescopic paths between $x$ and
$y$.
\end{definition}

\begin{remark}
Let $p$ be a $v$-telescopic path. The vertical length of each
maximal $v$-vertical segment in $p$ is greater or equal to $1$.

Any point in $\tilde{X}$ is at vertical distance smaller than
$\frac{1}{2}$ from a stratum over a vertex of $\mathcal T$. Thus,
when dealing with the behavior of (quasi)-geodesics or with the
hyperbolicity of $\tilde{X}$, there is no harm in requiring that
telescopic paths begin and end at strata over vertices of $\mathcal
T$, as was done in Definition \ref{ciel}.

For the sake of simplification, we will often forget the exponents
in the vertical, horizontal and telescopic lengths, unless some
ambiguity might exist.
\end{remark}

\begin{lemma}
\label{indispensable2} Let $(\tilde{X},{\mathcal T},\pi)$ be a
tree of hyperbolic spaces.

\begin{enumerate}
  \item There exist $\lambda_+(v) \geq 1$, $\mu(v) \geq 0$ such that,
if ${w}_{0}$ and ${w}_{1}$ are any two $v$-vertical segments, with
initial (resp. terminal) points $x_{0}, x_{1}$ (resp. $y_{0},
y_{1}$) and such that $\pi({w}_{0}) = \pi({w}_{1}) = [a, b]$ then:

  $$\frac{1}{\lambda^{d_{\mathcal T}(a,b)}_+(v)} d_{hor}(x_{0},x_{1}) - \mu(v) \leq
  d_{hor}(y_{0},y_{1}) \leq \lambda^{d_{\mathcal T}(a,b)}_+(v)
  d_{hor}(x_{0},x_{1}) + \mu(v)$$

The constants $\lambda_{+}(v), \mu(v)$ will be referred to as the
{\em constants of quasi-isometry}.

  \item  $\displaystyle \lim_{n \to + \infty} d_{hor}(x_0,x_n) = + \infty \Leftrightarrow
  \lim_{n \to + \infty} d^v_{tel}(x_0,x_n) = + \infty$ whenever ${(x_n)}_{n \in \mz^+}$
  is a sequence of points in some stratum.

  \item For any $v,v^\prime \geq 0$,
  $(\tilde{X},d^{v}_{tel})$ is quasi isometric to
  $(\tilde{X},d^{v^\prime}_{tel})$.

  \item For any $\alpha, \beta \in {\mathcal T}$ and $v
\geq 0$ there exists $C(v,d_{{\mathcal T}}(\alpha,\beta))$,
increasing in both variables, such that for any $x, y, z \in
X_{\alpha}$ with $z \in [x,y]$, whenever $x^\prime,y^\prime,z^\prime
\in X_\beta$ are the endpoints of $v$-vertical segments starting
respectively at $x,y$ and $z$, then $z^\prime \in
\vois{C(v,d_{{\mathcal T}}(\alpha,
\beta))}{hor}{[x^\prime,y^\prime]}$.
  \item For any $0 \leq w$, there is $D(w)$ such that, if $s$
  is a $v$-vertical segment, then $s$ is a $(D(w),D(w))$-quasi
  geodesic for the $w$-telescopic
  distance.
\end{enumerate}
\end{lemma}

\pols{Item (a) is a straightforward consequence of the definition of
a vertical segment. Items (b) and (c) are consequences of the
existence of the constants of quasi isometry given by the first
item. Item (d) amounts to saying that the image of a geodesic under
a $(a,b)$-quasi isometric embedding is $C(a, b)$-close to any
geodesic between the images of the endpoints. This is a well-known
assertion, see for instance \cite{CDP}. Like Item (a), Item (e) is
checked by a straightforward computation.}

\begin{remark}
Throughout all the text, the constants appearing in each lemma,
corollary or proposition will be denoted by $C,D,\cdots$ and
thereafter they will be referred to by the same letter with the
number of the lemma, corollary or proposition in subscript. For
instance, if Lemma 3.4 introduces the constants $C$ and $D$, for
referring afterwards to these constants, we will write $C_{3.4}$ and
$D_{3.4}$.
\end{remark}

\subsection{Exponential separation of vertical segments,
hallways-flare property}

\begin{definition} (compare \cite{BF})
\label{definition principale}

A tree of hyperbolic spaces $(\tilde{X},{\mathcal T},\pi)$ satisfies
the {\em hallways-flare property} if for any $v \geq 0$ there exist
positive integers $t_{0},M$ and a constant $\lambda
> 1$ such that, for any $\alpha \in {\mathcal T}$, for any two points $\beta, \gamma
 \in \partial B_{\alpha}(t_0)$ with $d_{{\mathcal T}}(\beta,\gamma)
 = 2 t_0$, any two $v$-vertical segments $s_0,s_1$ over
 $[\beta, \gamma]$ with $d_{hor}(s_0 \cap X_\alpha,s_1 \cap X_\alpha) \geq M$
 satisfy:

 $$\mathrm{max}(d_{hor}(s_0 \cap X_\beta,s_1 \cap X_\beta),d_{hor}(s_0 \cap X_\gamma,s_1 \cap X_\gamma)) \geq \lambda d_{hor}(s_0 \cap
X_\alpha,s_1 \cap X_\alpha)$$
\end{definition}

We will sometimes say that the $v$-vertical segments are {\em
$M$-exponentially separated}. The constants $\lambda,M,t_0$ will be
referred to as the {\em constants of hyperbolicity}. We now state a
very general lemma about these constants.

\begin{lemma}
\label{indispensable} Let $(\tilde{X},{\mathcal T})$ be a tree of
hyperbolic spaces satisfying the hallways-flare property.
\begin{enumerate}
  \item The constants of hyperbolicity and quasi isometry can be
  chosen arbitrarily large enough.
  \item If  $\lambda,M,t_0$ are the constants of hyperbolicity and $M$ is chosen sufficiently large enough,
  then there
exists $C$ such that, for any  $\alpha \in {\mathcal T}$,
  for any $\beta, \gamma \in \partial B_{t_0}(\alpha)$ with
  $\alpha \in [\beta,\gamma]$, for any two $v$-vertical
  segments $s_0,s_1$ over $[\beta,\gamma]$ such that
  $d_{hor}(x_0,x_1) \geq M$ where $x_i= s_i \cap X_\alpha$, if
  the endpoints $y_0,y_1$ of $s_0,s_1$ in $X_\beta$ (resp. in
  $X_\gamma$)
   satisfy: $$\frac{1}{\lambda} d_{hor}(x_0,x_1) < d_{hor}(y_0,y_1),$$ then, for
any $n \geq 1$, for any ${\mathcal T}$-geodesic $w$
starting at $\alpha$ with $[\alpha,\beta] \subset
  w$ (resp. $[\alpha,\gamma] \subset
  w$) and $|w|_{\mathcal T} \geq C + n t_0$: $$d^i_{hor}(wx,wy) \geq
  \lambda^n d^i_{hor}(x,y).$$
\end{enumerate}
\end{lemma}

The hallways-flare property above requires that the exponential
separation of the $v$-vertical segments be satisfied for {\em any}
$v \geq 0$. It suffices in fact that it be satisfied for {\em some}
$v$ sufficiently large enough as we are now going to check (see
Lemma \ref{bo2}).

\begin{definition}
\label{projection horizontale} Let $\tilde{X}$ be a tree of
hyperbolic spaces and let $S$ be a horizontal subset which is
(quasi) convex in its stratum, for the horizontal metric. If $x$ is
any point in $\tilde{X}$ then $P^{hor}_{S}(x)$ denotes any point $y$
in $S$ such that $d_{hor}(x,y) \leq d^{i}_{hor}(x,{S}) + 1$
\end{definition}

\begin{lemma}
\label{preparation} There exists $C$ such that if $v \geq C$, if $h$
is a horizontal geodesic in a stratum over some vertex $w$ of
$\mathcal T$, if $e$ is an edge of $\mathcal T$ incident to $w$ such
that no $v$-vertical segment starting at $h$ can be defined over
$e$, then $\mathrm{diam}_{X_w}(P^{hor}_h(i_{e,w}(X_e))) \leq 2
\delta$, where $\delta$ is the constant of hyperbolicity of the
strata and $i_{e,w}(X_e)$ denotes the quasi isometric embedding of
the edge-space $X_e$ into the vertex-space $X_w$.
\end{lemma}

\pols{The $\delta$-hyperbolicity of the strata for the horizontal
metric gives a constant $c$, depending on the constants of quasi
isometry, such that, for any two points $x,y \in i_{e,w}(X_e)$, any
horizontal geodesic $[x,y]$ lies in the horizontal $c$-neighborhood
of $i_{e,w}(X_e)$. Choose $v > 2 \delta + c$. Assuming that no
$v$-vertical segment starting at $h$ can be defined over $e$, since
horizontal geodesic rectangles are $2 \delta$-thin, we get $[x,y]
\cap \vois{2 \delta}{hor}{h} = \emptyset$ for any two points $x,y
\in i_{e,w}(X_e)$ and any horizontal geodesic $[x,y]$. The
conclusion follows by the $2 \delta$-thinness of the geodesic
rectangles.}

\begin{lemma}
\label{bo2} Let $(\tilde{X},{\mathcal T},\pi)$ be a tree of
hyperbolic spaces. If $v \geq C_{\ref{preparation}}$ is such that
the $v$-vertical segments of $\tilde{X}$ are exponentially separated
with constants of hyperbolicity $\lambda_v
> 1$, $M_v,t_0 \geq 0$ then for any $w \geq v$, the
$w$-vertical segments are exponentially separated, with constants of
hyperbolicity $\lambda_w > 1$, $M_w \geq 0$ and $t_0$.
\end{lemma}

\pols{Consider $\alpha, \beta, \gamma$ in $\mathcal T$ with $\alpha
\in [\beta,\gamma]$ and $d_{{\mathcal T}}(\alpha,\beta) =
d_{{\mathcal T}}(\alpha,\gamma) = t_0$. Consider two $w$-vertical
segments $S_0,S_1$ over
 $[\beta, \gamma]$ with $d_{hor}(x_0,x_1) \geq M$, where $x_i = S_i \cap
 X_\alpha$ and $M > M_v$. We distinguish two cases:

\begin{itemize}
  \item there exist $v$-vertical segments $s_0,s_1$ passing through
  $x_0,x_1$ and defined over $[\beta,\gamma]$. From Item (a) of
  Lemma \ref{indispensable2}, each endpoint of the $s_i$'s is at
  bounded horizontal distance from an endpoint of a $S_i$, where the
  upper-bound only depends on $w, t_0$ and the constants of quasi
  isometry. Thus choosing $M$ sufficiently large enough with respect
  to $w$ gives the desired inequality between $d_{hor}(x_0,x_1)$ and
  $\mathrm{max}(d_{hor}(S_0 \cap X_\beta,S_1 \cap X_\beta), d_{hor}(S_0 \cap
  X_\gamma,S_1 \cap X_\gamma))$.
  \item the other case: since $v$ has been chosen greater than
  $C_{\ref{preparation}}$, there is some stratum $X_\mu$, $\mu \in
  [\beta,\gamma]$ such that $d_{hor}(S_0 \cap X_\mu,S_1 \cap X_\mu)$
  is bounded above by a constant depending on $w,\delta,t_0$ and the
  constants on quasi isometry. By Item (a) of Lemma
  \ref{indispensable2}, we get an upper-bound on $d_{hor}(x_0,x_1)$.
  Setting $M$ greater that this upper-bound, we get the lemma.
  \end{itemize}}

\section{Approximation of quasi geodesics: a ``simple'' case}
\label{aqgpc}
From a group-theoretical point of view, the case
treated in this section allows one to deal with semi-direct products
of (relatively) hyperbolic groups with free groups but not with
HNN-extensions and amalgamated products along proper subgroups. For
this we need the similar, but more general, result of Section
\ref{aqggc}.

Beware that the corridors (and further the generalized and pseudo
corridors) defined below are not the hallways of \cite{BF}. The
reason is that we are interested in exhibiting quasi convex subsets
of our trees of hyperbolic spaces and the hallways of \cite{BF}, in
general, are not quasi convex.

\begin{definition}
Let $(\tilde{X},{\mathcal T},\pi)$ be a tree of hyperbolic spaces.
Let $\sigma_1,\sigma_2 \colon T \rightarrow
  \tilde{X}$ be two maximal (in the sense of the inclusion) $v$-vertical trees.
  A union of horizontal geodesics, at most one
  in each stratum, connecting each point of $\sigma_1(T)$ to a point of $\sigma_2(T)$
  is a {\em $v$-corridor}.
\end{definition}

\begin{remark}
Let $(\tilde{X},{\mathcal T},\pi)$ be a tree of hyperbolic spaces
the attaching-maps of which are all quasi isometries. Then, as soon
as $v \geq C_{\ref{preparation}}$, given any two points $x,y$ in
$\tilde{X}$, there is a
 $v$-corridor $\mathcal C$ whose
vertical boundaries pass through $x$ and $y$. Moreover
$\pi({\mathcal C}) = {\mathcal T}$.
\end{remark}

\begin{definition}
Let $\mathcal C$ be any subset of a tree of hyperbolic spaces
$\tilde{X}$ which is a union of horizontal geodesics, at most one in
each stratum (for instance $\mathcal C$ might be a corridor). Let
$X_\alpha$ be some stratum of $\tilde{X}$ and let $x$ be any point
in $X_\alpha$. The notation $P^{hor}_{{\mathcal C}}(x)$ stands for
$P^{hor}_{{\mathcal C} \cap X_\alpha}(x)$ (see Definition
\ref{projection horizontale}).
\end{definition}

Before stating Lemma \ref{casi leafe} below, we would like to insist
on two points:

\begin{itemize}
  \item The projection $P^{hor}_{\mathcal C}$ is a projection {\em in the
  strata} which only refers to the horizontal metric defined on each stratum.
  \item Item (b) does not tell anything about the behavior of the
telescopic (quasi)-geodesics in a tree of hyperbolic spaces. It only
allows one to consider a corridor as a quasi geodesic telescopic
metric space.
\end{itemize}

\begin{lemma}
\label{casi leafe} Let $(\tilde{X},{\mathcal T},\pi)$ be a tree of
hyperbolic spaces. Let $\mathcal C$ be a $v$-corridor (or a
generalized $v$-corridor - see Definition \ref{generalized
corridor}) in $\tilde{X}$ ($v \geq C_{\ref{preparation}}$). Then:
\begin{enumerate}
\item There exists $C(v) \geq v$ such that, if $s$ is a $v$-vertical
 segment, then $P^{hor}_{\mathcal C}(s)$ is a
$C(v)$-vertical segment.
\item For any $w \geq C(v)$, $({\mathcal C},d^{w}_{tel})$ is a quasi
geodesic metric space.
\end{enumerate}
\end{lemma}

\pols{If $\sigma \colon w \rightarrow \tilde{X}$ is the section of
$\pi$ such that $s = \sigma(w)$ then $P^{hor}_{\mathcal C}(s)$ is
the image of $w$ under the map $P^{hor}_{\mathcal C} \circ \sigma$.
This map is a section of $\pi$ since the projection
$P^{hor}_{\mathcal C}$ is a projection in each stratum. We want to
prove that $P^{hor}_{\mathcal C} \circ \sigma$ is a quasi isometric
embedding of $w$ into $(\tilde{X},d_{tel})$. Assume $w$ is a single
edge. Since $v \geq C_{\ref{preparation}}$ and since $\mathcal C$ is
a (generalized) $v$-corridor, if it is defined over $w$ then
$v$-vertical segments can be defined over $w$ starting at each point
of ${\mathcal C} \cap X_{i(w)}$. Let $\sigma_0 \colon w \rightarrow
\tilde{X}$ be such a $v$-vertical segement starting at $P_{{\mathcal
C}}(\sigma(i(w)))$. By Items (a) and (c) of Lemma
\ref{indispensable2}, $d_{hor}(P_{{\mathcal C}}(\sigma(t(w))),
\sigma_0(P_{{\mathcal C}}(\sigma(i(w))))$ is bounded above by a
constant. Thanks to Item (b) of Lemma \ref{indispensable2}, this
proves Item (a) of the current lemma. Item (b) is now easy.}

\begin{definition}
\label{def de diagonale} A {\em diagonal} is a horizontal geodesic
which minimizes the horizontal distance between two vertical trees
passing through its endpoints.
\end{definition}

\begin{theorem}
\label{pipeau} Let $\tilde{X}$ be a tree of hyperbolic spaces which
satisfies the hallways-flare property. Assume that each
attaching-map from an edge-space into a vertex-space is a quasi
isometry.

Choose a constant $L > 0$ greater than some critical constant. Then
there are $C(L,a,b),$ $D(L,a,b) \geq 0$ such that for any telescopic
$(a,b)$-quasi geodesic $g$ in $\tilde{X}$ there is a telescopic
$(D(L,a,b),$ $D(L,a,b))$-quasi geodesic ${\mathcal G}$ satisfying the
following properties:

\begin{enumerate}
  \item $d^H_{tel}(g,{\mathcal G}) \leq C(L,a,b)$;
  \item ${\mathcal G}$ is contained in a corridor
  $\mathcal C$ the vertical boundary trees of which pass through
  the endpoints of $g$;
  \item at the exception of at most one, each maximal horizontal subpath of ${\mathcal G}$
  is a diagonal with horizontal length
greater or equal to $L$ whereas the last maximal horizontal subpath
has horizontal length less or equal to $L$;
  \item the corridor $\mathcal C$ only depends on the endpoints of
  $g$;
  \item at the exception
  of its first and last maximal vertical segments, which depend on the initial
  and terminal points of $g$, ${\mathcal G}$ only
  depends on the choice of the corridor $\mathcal C$.
\end{enumerate}

\end{theorem}

\pot{pipeau}{We need two important propositions, which we state now
but the proofs of which are postponed for a while.

\begin{proposition}
\label{proposition importante 1} Let ${\mathcal C}$ be a corridor
(or a generalized corridor - see Definition \ref{generalized
corridor}) in a tree of hyperbolic spaces. Assume that $\mathcal C$
satisfies the hallways-flare property. Then there exists $C(L,a,b)$
such that, if $L$ is the horizontal distance in some stratum
$X_\alpha$ between two $v$-vertical trees $T_1,T_2$, if $\mathcal G$
is a $v$-telescopic $(a,b)$-quasi geodesic of $({\mathcal
C},d^{v}_{tel})$ from $T_1$ to $T_2$ which starts or ends at a
stratum where the horizontal distance between $T_1$ and $T_2$ is
greater than $L$, then $\mathcal G$ is contained in the
$C(L,a,b)$-neighborhood of the union of the vertical segments
connecting its endpoints to the points $T_i \cap X_\alpha$. The
constant $C(L,a,b)$ is increasing with $L$ as soon as $L$ is greater
than some critical constant.
\end{proposition}

See Section \ref{propimpo1} for a proof.

\begin{proposition}
\label{quasiconvexite} Let $\tilde{X}$ be a tree of hyperbolic
spaces which satisfies the hallways-flare property and the
attaching-maps of which are quasi isometries. There exists $C(a,b)$
such that, if $g$ is a telescopic $(a,b)$-quasi geodesic in
$\tilde{X}$, if $\mathcal C$ is a corridor containing the endpoints
of $g$ then
$$g \subset \vois{C(a, b)}{tel}{{\mathcal C}}.$$
\end{proposition}

See Section \ref{quasi convexe} for a proof.

\begin{lemma}
\label{shot} Let $\mathcal C$ be a (generalized) corridor in a tree
of hyperbolic spaces. There exists $C \geq 0$ such that, for any two
points $x,y$ in a same stratum ${X}$, $d_{hor}(P^{hor}_{{\mathcal
C}}(x),P^{hor}_{{\mathcal C}}(y)) \leq d_{hor}(x,y) + C$. The same
inequality holds when projecting $x,y$ to the image of the embedding
of an edge-space into a vertex-space.
\end{lemma}

\pols{Since strata are $\delta$-hyperbolic space for the horizontal
metric and the subspaces to which one projects are (quasi) convex
subsets of their stratum for this horizontal metric, this is a
consequence of \cite{CDP}, Corollary 2.2 page 109.}

\begin{lemma}
\label{dimanche} Let $g$ be a $v$-telescopic path, which is a
$(a,b)$-quasi geodesic for the telescopic distance $d^{v}_{tel}$.
Let $\mathcal C$ be a (generalized) corridor. Then there exists
$C(a,b,r) \geq 1$ such that, if $g \subset \vois{r}{hor}{{\mathcal
C}}$ then $P^{hor}_{\mathcal C}(g)$ is a $C_{\ref{casi
leafe}}(v)$-telescopic $(C(a,b,r),C(a,b,r))$-quasi geodesic of
$({\mathcal C},d^{C_{\ref{casi leafe}}(v)}_{tel})$.
\end{lemma}

\pols{Lemma \ref{casi leafe} implies in a straightforward way that
$P^{hor}_{\mathcal C}(g)$ is a $C_{\ref{casi leafe}}(v)$-telescopic
path. Let us consider any two points $x,y$ in ${\mathcal G} =
P^{hor}_{\mathcal C}(g)$. There are $r$-close to two points
$x^\prime,y^\prime$ in $g$. We denote by $g_{x^\prime y^\prime}$ the
subpath of $g$ between these last two points and by ${\mathcal
G}_{xy}$ the subpath of $\mathcal G$ between $x$ and $y$. Since we
now consider the $C_{\ref{casi leafe}}(v)$-telescopic distance,
$|{\mathcal G}_{xy}|^{C_{\ref{casi leafe}}(v)}_{vert} = |g_{x^\prime
y^\prime}|^v_{vert}$. From Lemma \ref{shot} and since any two
maximal horizontal subpaths of $\mathcal G$ are separated by a
vertical segment of vertical length at least $1$, we then get
$|{\mathcal G}_{xy}|^{C_{\ref{casi leafe}}(v)}_{tel} \leq 2
C_{\ref{shot}} |g_{x^\prime y^\prime}|^{v}_{tel}$. Since $g$ is a
$v$-telescopic $(a,b)$-quasi geodesic, $|g_{x^\prime
y^\prime}|^{v}_{tel} \leq a d^{v}_{tel}(x^\prime,y^\prime) + b$. But
$d^{v}_{tel}(x^\prime,y^\prime) \leq 2r + d^{v}_{tel}(x,y)$.
Therefore: $$|{\mathcal G}_{xy}|^{C_{\ref{casi leafe}}(v)}_{tel}
\leq 2 C_{\ref{shot}} (a(2r + d^{v}_{tel}(x,y))+b).$$ Since all
telescopic distances are quasi isometric (Item (c) of Lemma
\ref{indispensable2}), we so get the right inequality for the quasi
geodesicity of
$P^{hor}_{\mathcal C}(g)$. We leave the reader work out the similar proof of the left inequality.} \\

Let $(\tilde{X},{\mathcal T},\pi)$ be a tree of hyperbolic spaces.
Choose $v \geq C_{\ref{preparation}}$. Let $g$ be a $v$-telescopic
$(a,b)$-quasi geodesic. Since the attaching maps of the tree of
hyperbolic spaces are all quasi isometries, there is a $v$-corridor
$\mathcal C$ the vertical boundaries of which pass through the
endpoints of $g$. This corridor $\mathcal C$ satisfies
$\pi({\mathcal C}) = {\mathcal T}$. From Proposition
\ref{quasiconvexite}, $g \subset \vois{C_{\ref{quasiconvexite}}(a,
b)}{tel}{{\mathcal C}}$. Since $\pi({\mathcal C}) = {\mathcal T}$,
Item (b) of Lemma \ref{indispensable2} then implies the existence of
$C^\prime(a,b)$ such that $g \subset
\vois{C^\prime(a,b)}{hor}{{\mathcal C}}$. From Lemma \ref{dimanche},
${\mathcal G} \equiv P^{hor}_{\mathcal C}(g)$ is a $C_{\ref{casi
leafe}}(v)$-telescopic $(A,A)$-quasi geodesic of $({\mathcal
C},d^{C_{\ref{casi leafe}}(v)}_{tel})$, with $A \equiv
C_{\ref{dimanche}}(a,b,C^\prime(a, b))$. From Lemma \ref{bo2},
$\mathcal C$ satisfies the hallways-flare property, more precisely
the $C_{\ref{casi leafe}}(v)$-vertical segments are exponentially
separated. From Item (b) of Lemma \ref{indispensable}, this implies
in particular that the endpoints of any diagonal with horizontal
length greater than some constant $M$ are exponentially separated in
all the directions of $\mathcal T$ outside a region with vertical
width bounded by $2 C_{\ref{indispensable}}$.

Choose $L \geq M$. Consider a diagonal $h_0$ with horizontal length
$L$ from a vertical boundary of $\mathcal C$ to a $C_{\ref{casi
leafe}}(v)$-vertical tree in $\mathcal C$. The quasi geodesic
$\mathcal G$ joins these two vertical trees of $\mathcal C$, let
${\mathcal G}_0$ be the corresponding subpath of $g$. From
Proposition \ref{proposition importante 1}, ${\mathcal G}_0$ is
contained in the $C_{\ref{proposition importante
1}}(L,A,A)$-neighborhood of the union of the vertical segments
$s_0,s_1$ from the endpoints of ${\mathcal G}_0$ to those of $h_0$.
From our observation above about the exponential separation of the
endpoints of $h_0$, there is some $\kappa > 0$ such that, outside
the region in $\mathcal C$ centered at $h_0$ with vertical width
$\kappa$, the horizontal geodesics between the vertical trees of the
endpoints of $h_0$ have horizontal length greater than $3
C_{\ref{proposition importante 1}}(L,A,A)$. We so get a constant $K
> 0$, not depending of the quasi geodesic considered, such that
$d^H_{tel}({\mathcal G}_0,s_0 \cup h_0 \cup s_1) \leq K$.

By continuing the construction of diagonals $h_1,\cdots,h_r$ as was
constructed $h_0$, at each step starting from the last vertical tree
considered, we eventually get an approximation of a maximal subpath
${\mathcal G}^{\prime}$ of ${\mathcal G}$ by a concatenation of
diagonals and vertical segments between these diagonals as was
announced by Theorem \ref{pipeau}. Observe that:

\begin{itemize}
  \item The corridor $\mathcal C$ only depends on the choices made
  for the vertical trees through the endpoints of $g$.
  \item The diagonals only depend on $\mathcal C$, and not on $g$.
\end{itemize}

The subpath ${\mathcal G}^\prime$ above is characterized by the fact
that there is no diagonal with horizontal length $L$ between any
vertical tree in $\mathcal C$ through its terminal point and the
vertical boundary containing the terminal point of $g$. The choice
of the last horizontal geodesic then only depends on the vertical
position of the endpoint of the last diagonal $h_k$: if the terminal
point of $h_k$ is at horizontal distance smaller than $L$ of the
terminal vertical boundary, then choose this horizontal geodesic as last
one; otherwise go along a vertical segment to the nearest stratum
where the horizontal distance between the two considered vertical
trees is equal to $L$.

The path we so get satisfies all the properties announced by Theorem
\ref{pipeau}.}

\section{Approximation of quasi geodesics: the general case}

\label{aqggc}

In order to give a simple statement, we added in Theorem
\ref{pipeau} the restriction that the attaching-maps of the tree of
spaces be quasi isometries, instead of quasi isometric embeddings.
We now come to the more general statement.

 \begin{definition}
 \label{generalized corridor}
 A {\em generalized $v$-corridor $\mathcal C$}
  is a union of horizontal geodesics, at most one in each stratum, such
  that $\pi({\mathcal C}) \equiv T$ is a subtree of $\mathcal T$ which admits
  a decomposition in subtrees $T_i$, with $T_i \cap T_j$ either empty or reduced to a single
  point when $i \neq j$, satisfying the following properties:

  \begin{enumerate}
    \item for each $i$, ${\mathcal C} \cap \pi^{-1}(T_i) \equiv {\mathcal
    C}_i$ is a union of horizontal geodesics between two vertical
    trees,
    \item if $v$ is a vertex of $\mathcal T$ in $T$ and
    $e$ is an edge of $\mathcal T$ which is incident to $v$ but does not belong to
    $T$, then there is no $v$-vertical segment over $e$ starting from $\mathcal
    C$;
    \item if $x$ is a point in the horizontal boundary of some
    ${\mathcal C}_i$ such that some $v$-vertical segment $s$ with $\pi(s) \subset
    T_j$, $j \neq i$, starts
    at $x$, then $x$ is in ${\mathcal C}_j$.
  \end{enumerate}
\end{definition}

\begin{definition}
\label{pseudo corridor}
A {\em pseudo-corridor} is a concatenation of
generalized corridors ${\mathcal C}_i$ and of horizontal geodesics
$h_j$ such that:
  \begin{itemize}
    \item either $\pi({\mathcal C}_i) \cap \pi({\mathcal C}_k)$ is
    reduced to a single point and then there is exactly one $h_j$
    connecting ${\mathcal C}_i$ to ${\mathcal C}_k$,
    \item or $\pi({\mathcal C}_i) \cap \pi({\mathcal C}_k)$ is
    empty.
  \end{itemize}
\end{definition}

\begin{remark}
If $\tilde{X}$ is a tree of hyperbolic spaces
then, if  $v \geq C_{\ref{preparation}}$,
given any two points $x,y$ in $\tilde{X}$, there is a
pseudo-corridor $\mathcal C$ whose vertical boundaries pass through
$x$ and $y$.
\end{remark}

\begin{theorem}
\label{pipeau50} Let $\tilde{X}$ be a tree of hyperbolic spaces
which satisfies the hallways-flare property. The conclusions of
Theorem \ref{pipeau} remain true, when dropping the extra-hypothesis
on the attaching-maps, with the following modifications:

\begin{itemize}
  \item one substitutes the word ``corridor'' by the word
``pseudo-corridor'',
  \item the maximal horizontal subpaths of $\mathcal G$ are
  diagonals with horizontal length greater than $L$ at the exception
  of at most one in each generalized corridor forming the
  pseudo-corridor.
\end{itemize}
\end{theorem}

\pot{pipeau50}{We first need an adaptation of Proposition
\ref{quasiconvexite} to this more general setting:

\begin{proposition}
\label{qc}
Proposition \ref{quasiconvexite} remains true when
dropping the assumption on the attaching-maps if one substitutes the
word ``corridor'' by the word ``generalized corridor''.
\end{proposition}

See the proof in subsection \ref{popqc} of Section \ref{quasi
convexe}. Unfortunately, two distinct points do not necessarily
belong to a same generalized corridor. This is why we needed to
introduce the pseudo-corridors, and why we need Lemma \ref{testons2}
below.

\begin{lemma}
\label{testons2} Let $\tilde{X}$ be a tree of
hyperbolic spaces. There is a $C(a,b) \geq 0$ such that, if $g$ is
any telescopic $(a,b)$-quasi geodesic, then $g$ lies in the
telescopic $C(a,b)$-neighborhood of a pseudo-corridor. This
pseudo-corridor only depends on the endpoints of $g$.
\end{lemma}

\pols{Let $X_\alpha,X_\beta$ be the strata containing the initial
and terminal points of $g$. There is a unique sequence of $\gamma_i
\in [\alpha,\beta]$, $i=1,\cdots,k \mbox{, } \gamma_0 = \alpha,
\gamma_k = \beta$ such that $\gamma_i$ maximizes $d_{\mathcal
T}(\gamma_{i-1},\phi)$ among all $\phi$'s in $[\alpha,\beta]$ such
that some $v$-vertical segment connects $X_{\gamma_{i-1}}$ to
$X_\phi$. We denote by $Y_i$ the maximal region of $X_{\gamma_i}$
for which $v$-vertical segments are defined from $Y_i$ to
$X_{\gamma_{i+1}}$. Since $v \geq C_{\ref{preparation}}$, from Lemma
\ref{preparation}, for any $\gamma_i$, $[\gamma_{i} \gamma_{i+1}]
Y_i$ is connected to $Y_{i+1}$ by a horizontal rectangle $R_i$ of
width at most $2 \delta$. We denote by $a_i$ (resp. $b_i$) a point
in $R_i \cap [\gamma_{i} \gamma_{i+1}] Y_i$ (resp. in $R_i \cap
Y_{i+1}$). We then denote by ${\mathcal C}_i$ a generalized
$v$-corridor between $b_i$ and $a_{i+1}$ (that is containing a
horizontal geodesic between $b_i$ and the intersection-point of
$X_{\gamma_{i+1}}$ with a vertical tree through $a_{i+1}$). We set
$h_i = [a_i,b_i]$. The $(a,b)$-quasi geodesic $g$ connects a point
in the horizontal $2 \delta$-neighborhood of $b_i$ to a point in the
horizontal  $2 \delta$-neighborhood of $a_{i+1}$. Let us denote by
$g_i$ such a subpath of $g$. Thus, by connecting the endpoints of
$g_i$ to $b_i$ and $a_{i+1}$ we obtain a $(a,b+2 \delta)$-quasi
geodesic between $b_i$ and $a_{i+1}$, still denoted by $g_i$. From
Proposition \ref{qc}, ${\mathcal G}_i$ lies in the telescopic
$C_{\ref{qc}}(a,b+4 \delta)$-neighborhood of ${\mathcal C}_i$.
Obviously, since the width of $R_i$ is less or equal to $2 \delta$,
the subpath of $g$ between two $g_i$'s is in the $2
\delta$-neighborhood of $h_i$. This completes the proof of Lemma
\ref{testons2}.} \\

In order to follow the proof of Theorem \ref{pipeau} in the more
general setting we are confronted to, we still need an additional
result. Lemma \ref{testons} below allows one to substitute the given
quasi geodesic $g$ by a quasi geodesic $\mathcal G$ with the
following properties:

\begin{itemize}
  \item it has the same endpoints, and is
Hausdorff-close to $g$ with respect to the telescopic distance,
  \item it admits a decomposition in subpaths ${\mathcal G}_i$ such
  that both endpoints of ${\mathcal G}_i$ lie in a same generalized
  corridor ${\mathcal C}_i$ of the considered pseudo-corridor and
  $\pi({\mathcal G}_i) \subset \pi({\mathcal C}_i)$.
\end{itemize}

This last property is needed in order to apply Lemma \ref{dimanche}.

\begin{lemma}
\label{testons} Let $\tilde{X}$ be a tree of
hyperbolic spaces. If $v \geq C_{\ref{preparation}}$ then there
exists $C(a,b)$ such that, if $g$ is any telescopic $(a,b)$-quasi
geodesic the endpoints of which lie in a generalized $v$-corridor
$\mathcal C$, then there is a telescopic $(a,b + 2 \delta)$-quasi
geodesic $\mathcal G$ with $d^H_{tel}(g,{\mathcal G}) \leq C(a,b)$
and $\pi({\mathcal G}) \subset \pi({\mathcal C})$.
\end{lemma}

\pols{Let $\gamma \in {\mathcal T}$ be an endpoint of $\pi({\mathcal
C})$. Assume that $g^\prime$ is a maximal subpath
 of $g$ with endpoints in $X_\gamma$ and such that
$\pi(g^\prime) \cap \pi({\mathcal C}) = \gamma$. Then, since $v \geq
C_{\ref{preparation}}$, Lemma \ref{preparation} tells us that the
endpoints of $g^\prime$ are $2 \delta$-close with respect to the
horizontal distance. Since $g$ is a $(a,b)$-quasi geodesic,
$g^\prime$ is $(2 a \delta + b)$-close to $X_\gamma$ with respect to
the telescopic distance. Substituting $g^\prime$ by a horizontal
geodesic between its endpoints and repeating this substitution for
all the subpaths of $g$ as $g^\prime$ yields a quasi geodesic as
announced.} \\

With the above adaptations in mind, the proof of Theorem
\ref{pipeau50} is now a duplicate of the proof of Theorem
\ref{pipeau}.}

\section{Hyperbolicity and Weak Relative hyperbolicity}

\label{rh}

\subsection{Hyperbolicity of trees of spaces}

\label{hts}

Theorem \ref{main theorem one} generalizes Bestvina-Feighn's
combination to non-proper hyperbolic spaces. Bowditch proposed such
a generalization in \cite{Bo2}.

\begin{theorem}
\label{main theorem one} Let $\tilde{X}$ be a tree of hyperbolic
spaces which satisfies the hallways-flare property. Then $\tilde{X}$
is a Gromov-hyperbolic metric space.
\end{theorem}

\pot{main theorem one}{We begin by proving the

\begin{theorem}
\label{theorem one bis} Let $\tilde{X}$ be a tree of
hyperbolic spaces which satisfies the hallways-flare property. There
exists $C(a,b)$ such that telescopic $(a,b)$-quasi geodesic bigons
are $C(a,b)$-thin.
\end{theorem}

\pots{Let  $g_{0}, g_{1}$ be the two sides of a telescopic
 $(a, b)$-quasi geodesic bigon. By Theorem \ref{pipeau50} (Theorem \ref{pipeau} suffices in
the case where attaching-maps of $\tilde{X}$ are quasi isometries),
there is a telescopic path
 $\mathcal G$ such that for $i=0,1$ and $r$ a constant chosen sufficiently large
enough, we have $d^{H}_{tel}(g_{i}, {\mathcal G}) \leq
 C_{\ref{pipeau50}}(r,a,b)$. Hence $d^{H}_{tel}(g_{0}, g_{1}) \leq 2
 C_{\ref{pipeau50}}(r,a,b)$
and Theorem \ref{theorem one bis}
is proved.} \\

The following lemma was first indicated
to the author by I. Kapovich:

\begin{lemma} \cite{Ga}
\label{ilya gautero} Let $(X,d)$ be a $(r,s)$-quasi geodesic space.
If for any $r^{\prime} \geq r, s^{\prime}
 \geq s$, there exists
 $\delta(r^{\prime},s^{\prime})$,
 such that  $(r^{\prime}, s^{\prime})$-quasi geodesic bigons are
 $\delta(r^{\prime}, s^{\prime})$-thin, then
$(X, d)$ is a $2\delta(r,3s)$-hyperbolic space.
\end{lemma}

 Theorem
\ref{theorem one bis} together with Lemma \ref{ilya gautero}
imply Theorem \ref{main theorem one}.}

\subsection{Weak relative hyperbolicity}

We first recall the definition. If $S$ is a discrete set, the {\em
cone with base $S$} is the space $S \times [0,\frac{1}{2}]$ with $S
\times \{0\}$ collapsed to a point, the {\em vertex of the cone}.
This cone is considered as a metric space, with distance function
$d_S((x,t),(y,t^\prime)) = t + t^\prime$. Let $(X, d)$ be a
quasigeodesic space. Putting a cone over a discrete subset
 $S$ of $X$ consists of pasting to $X$ a cone with base $S$ by
 identifying $S \times \{1/2\}$ with $S \subset X$. The resulting metric
 space, called the {\em coned space},
 $(\hat{X}_S,{d}_S)$ is such that all the points in $S$ are now at distance
 $\frac{1}{2}$ from the vertex of the cone and so at distance $1$ one
 from each other. The metric of the coned space is the {\em coned}, or
 {\em relative}, {\em metric}. If $\mathcal S$ is a
 disjoint union of sets, then the coned space  $X_{{\mathcal S}}$
 is the space obtained by putting a cone over each set in  $\mathcal S$.

\begin{definition} \cite{Fa}
A quasi geodesic space $(X,d)$ is {\em weakly hyperbolic relative to
a family of subsets $\mathcal S$} if the coned space
$(\hat{X}_{\mathcal S},d_{\mathcal S})$ is Gromov hyperbolic.

Let $G$ be a group with finite generating set $S$ and associated
Cayley graph $\Gamma_{G}$, and let ${\mathcal H} = \{H_1,\cdots \}$
be a (possibly infinite) family of infinite subgroups $H_i$ of $G$.

The group  $G$ is {\em
  weakly hyperbolic relative to $\mathcal H$} if $\Gamma_{G}$ is weakly hyperbolic
relative to the family of the right classes $x H_i$.

The subgroups $H_i$ in the family $\mathcal H$ are the {\em
parabolic subgroups} of $G$.
\end{definition}

\begin{remark}
The definition of Farb relative hyperbolicity given above is the
original one \cite{Fa}. It is equivalent to require that $G$
equipped with the relative metric,  i. e. the metric associated to
the system of generators $S \cup {\mathcal H}$, be hyperbolic.
However the introduction of the cones and of the coned Cayley graphs
above is needed to introduce farther in the paper the Bounded Coset
Penetration property.
\end{remark}

\begin{definition}
A {\em graph of weakly} {\em relatively hyperbolic groups} is a
graph of groups $({\mathcal G},{\mathcal H}_v,{\mathcal H}_e)$ such
that:
\begin{enumerate}
  \item Each edge group $G_e$ and each vertex group $G_v$ is
  weakly hyperbolic relative
  to a specified (possibly empty) finite family of
  infinite subgroups ${\mathcal H}_e$ and ${\mathcal H}_v$.
  \item For any edge $e$, $(G_e,|.|_{{\mathcal H}_e})$ is quasi
  isometrically embedded in $(G_{i(e)},|.|_{{\mathcal
  H}_{i(e)}})$ and in $(G_{t(e)},|.|_{{\mathcal
  H}_{t(e)}})$.
\end{enumerate}
\end{definition}

To a graph of groups $\mathcal G$ with fundamental group $\mathcal
J$, we associate a graph of spaces as follows:

\begin{itemize}
  \item each edge (resp. vertex) group $G_e$ (resp. $G_v$) is the
fundamental group of a {\em standard $2$-complex} $K_e$ (resp.
$K_v$): its $1$-skeleton is a rose, the petals of which are in
bijection with the generators of the group; its $2$-cells are glued
along the petals by simplicial maps of their boundaries, which
represent the relations;
  \item each edge space $K_e$ is glued to the vertex
spaces $K_{i(e)}$ and $K_{t(e)}$ by simplicial maps $\psi_{e,i(e)}$,
$\psi_{e,t(e)}$ which induce, on the level of the fundamental
groups, the injections of $G_e$ into $G_{i(e)}$ and into $G_{t(e)}$
coming with $\mathcal G$.
\end{itemize}

Let us consider the universal covering of this graph of spaces. This
is a tree of metric spaces $\pi \colon \tilde{X} \rightarrow
{\mathcal T}$ as
 defined in \ref{coqi}. The vertex (resp. edge) spaces are the universal covering of
the
 $K_v$'s (resp. $K_e$'s), these are just Cayley complexes for the
 edge and vertex groups of $\mathcal G$. They are equipped with the usual simplicial
 metric.

 In the case where $\mathcal G$ is
 a graph of weakly relatively hyperbolic groups, the edge and vertex
 groups
are weakly hyperbolic relative
 to certain subgroups. Associated to these subgroups is a relative metric.
 We equip the strata of the above constructed tree of spaces
 $\tilde{X}$
(the $1$-skeleton of a stratum is
 the Cayley graph of the corresponding edge or vertex group) with
 these coned metrics.
 We denote by $\widehat{X}$ the space obtained. This is a tree of
 hyperbolic spaces.

\begin{definition}
\label{inutile1} Let ${\mathcal G}$ be a graph of weakly relatively hyperbolic groups. The
universal covering of $\mathcal G$ satisfies the {\em relative
hallways-flare property} if the space $\widehat{X}$ constructed
above satisfies the hallways-flare property.
\end{definition}

\begin{theorem}
\label{intro sup} Let ${\mathcal G}$
be a finite graph of weakly relatively hyperbolic groups. If the
universal covering of $\mathcal G$ satisfies the relative
hallways-flare property, then the fundamental group of $\mathcal G$
is weakly hyperbolic relative to the family formed by all the
parabolic subgroups of the edge and vertex groups.
\end{theorem}

\pots{By Theorem \ref{main theorem one}, this readily follows from the definitions.} \\

The relatively hyperbolic automorphims we define below first
appeared in \cite{Ga1} where we announced a (weak) version of the
results of the present paper. They generalize the Gromov hyperbolic
automorphisms \cite{BF}.

\begin{definition}
\label{premieres definitions} Let $G = \langle S \rangle$ be a
finitely generated group and let ${\mathcal H} = \{H_1,\cdots,H_k\}$
be a finite family of subgroups of $G$.

\begin{enumerate}
  \item An automorphism $\alpha$ of $G$ is a {\em relative
automorphism of $(G,{\mathcal H})$} if $\mathcal H$ is
$\alpha$-invariant up to conjugacy, that is there is a permutation
$\sigma$ of $\{1,\cdots,k\}$ such that for any $H_{i} \in {\mathcal
H}$ there is $g_{i} \in G$ with $\alpha(H_{i}) = g^{-1}_{i}
H_{\sigma(i)}
 g_{i}$.
  \item The {\em $\mathcal
H$-word metric} $|.|_{\mathcal H}$ is the word-metric for $G$
equipped with the (usually infinite) set of generators which is the
union of $S$ with the elements of $G$ in the subgroups of the
collection $\mathcal H$.
   \item An automorphism $\alpha$ of $G$ is {\em hyperbolic relative to $\mathcal H$} if
   $\alpha$ is
a relative automorphism of $(G,{\mathcal H})$ and there exist
$\lambda
> 1$ and  $M,N \geq 1$ such that for any $w \in G$ with $|w|_{{\mathcal H}} \geq M$: $$\lambda
|w|_{\mathcal H} \leq \mathrm{max}(|\alpha^{N}(w)|_{\mathcal
H},|\alpha^{-N}(w)|_{\mathcal H}).$$
\end{enumerate}
\end{definition}

The definition of relatively hyperbolic automorphism given above is
slightly more general than the definition given in \cite{Ga1}. The
constant $M$ did not appear there. It is however more natural:
thanks to this additional constant $M$, the definition is obviously
invariant under conjugacy\footnote{The author is grateful to F.
Dahmani,
  V. Guirardel
and M. Lustig for this observation.}.

\begin{definition}
\label{auto relativement hyperbolique} Let $G$ be a finitely
generated group and let $\mathcal H$ be a finite family of subgroups
of $G$. A {\em uniform free group of relatively hyperbolic automorphisms}
\footnote{The author would like to thank M. Heusener for inciting
him to correct a previous formulation of this definition, which was
unnecessarily more restrictive.} of $(G,{\mathcal H})$ is a rank $r$
free group $\F{r}$ of relative automorphisms of $(G,{\mathcal H})$
such that there exist, for some (and hence any) basis $\mathcal A$
of $\F{r}$, $\lambda
> 1$ and $M,N \geq 1$ such that,
for any element $w \in G$ with $|w|_{{\mathcal H}} \geq M$, any pair
of automorphisms $\alpha, \beta$ with $|\alpha|_{{\mathcal A}} =
|\beta|_{{\mathcal A}} = N$ and $d_{{\mathcal A}}(\alpha,\beta) = 2
N$ satisfies:

$$\lambda |w|_{\mathcal H} \leq \mathrm{max}(|\alpha(w)|_{\mathcal H},|\beta(w)|_{\mathcal H}).$$
\end{definition}

\begin{theorem}
\label{theoreme principal faible} Let $G$ be a finitely generated
group and let $\mathcal H$ be
 a finite family of infinite subgroups  of $G$. Let $\F{r}$ be a
uniform free group of relatively hyperbolic automorphisms of $(G,
{\mathcal H})$. If $G$ is weakly hyperbolic relative to $\mathcal H$
then $G \rtimes \F{r}$
 is weakly hyperbolic relative to $\mathcal
 H$.
\end{theorem}

\pots{The group $G \rtimes \F{r}$ is the fundamental group of the
graph of groups which has $G$ as unique vertex group $G_v$, $G$ as
the $r$ edge-groups $G_{e_i}$ (the $e_i$'s are loops incident to
$v$) and the attaching endomorphisms of $G_{e_i}$ to $G_v$ are the
identity on one side and the automorphism $\alpha_i$ on the other
side, where the $\alpha_i$'s generate $\F{r}$. Since the
$\alpha_i$'s are relative automorphisms of $(G,{\mathcal H})$, each
one induces a quasi isometry from $(G_{e_i},{\mathcal H})$ to
$(G_v,{\mathcal H})$. Since $\F{r}$ is a uniform free group of
relatively hyperbolic automorphisms, the universal covering of this
graph of groups satisfies the relative hallways-flare property.
Theorem \ref{theoreme principal faible} is then a corollary of
Theorem \ref{intro sup}.} \\

From \cite{Ge2}, a hyperbolic group is weakly hyperbolic relative to
any finite family of quasi convex subgroups. We so get:

\begin{corollary}
\label{premier corollaire faible} Let $G$ be a hyperbolic group, let
$\mathcal H$ be a finite family of infinite subgroups of $G$ and let
$\alpha$ be an automorphism of $G$ which is hyperbolic relative to
$\mathcal H$. If $\mathcal H$ is quasi convex in $G$ then the
mapping-torus group $G_{\alpha} = G \rtimes_\alpha \mz$ is weakly
hyperbolic relative to $\mathcal H$.
\end{corollary}

\section{Strong relative hyperbolicity}

\label{strong}

Let $(\hat{X}_{\mathcal S},d_{\mathcal S})$ be a coned space (see
the beginning of the previous section) and let $\hat{g}$ be a
$(u,v)$-quasi geodesic in $(\hat{X}_{\mathcal S},d_{\mathcal S})$. A
{\em trace} $g$ of $\hat{g}$ in $(X,d)$ is obtained by substituting
each subpath of $\hat{g}$ not in $(X,d)$ by a subpath of $(X,d)$ in
$S$ with same endpoints, which is a geodesic for the metric induced
by $X$ on $S$. We say that  $g$ (or $\hat{g}$) {\em backtracks} if
$g$ reenters a subset in $\mathcal S$ that it left before.

\begin{definition}\cite{Fa}
A coned space $(\hat{X}_{\mathcal S},d_{\mathcal S})$ satisfies the
{\em Bounded-Coset Penetration property (BCP)} if there exists
$C(u,v)$ such that, for any two $(u,v)$-quasi geodesics
$\hat{g}_0,\hat{g}_1$ of $(\hat{X}_{\mathcal S},d_{\mathcal S})$
with traces $g_0,g_1$ in $(X,d)$, which have the same initial point,
which have terminal points at most $1$-apart and which do not
backtrack, the following two properties are satisfied:
\begin{enumerate}
  \item if both $g_0$ and $g_1$ intersects a set $S_i \in {\mathcal S}$ then their
first intersection points with $S_i$ are $C(u,v)$-close in $(X,d)$,
  \item if $g_0$ intersects a set $S_i$ that $g_1$ does not, then the length
in $(X,d)$ of $g_0 \cap S_i$ is smaller than $C(u,v)$.
\end{enumerate}
\end{definition}

\begin{definition}\cite{Fa}
A quasi geodesic space $(X,d)$ is {\em strongly hyperbolic relative
to a family of subsets $\mathcal S$} if the coned space
$(\hat{X}_{\mathcal S},d_{\mathcal S})$ is Gromov hyperbolic and
satisfies the BCP.

Let $G$ be a group with finite generating set $S$ and associated
Cayley graph $\Gamma_{G}$, and let ${\mathcal H} = \{H_1,\cdots, \}$
be a (possibly infinite) family of infinite subgroups $H_i$ of $G$.

The group  $G$ is {\em strongly hyperbolic relative to $\mathcal H$}
if $\Gamma_{G}$ is strongly hyperbolic relative to the union of the
right classes $x H_i$.
\end{definition}

\begin{definition}
\label{gsrh}
A {\em graph of strongly relatively hyperbolic groups}
is a graph of groups $({\mathcal G},{\mathcal H}_v,{\mathcal H}_e)$
such that:
\begin{enumerate}
  \item Each edge group $G_e$ and each vertex group $G_v$ is
  strongly hyperbolic relative
  to a specified (possibly empty) finite family of
  infinite subgroups ${\mathcal H}_e$ and ${\mathcal H}_v$.
  \item the edge collections
  ${\mathcal H}_e$ are required to be (possibly empty) families of conjugates of the
  subgroups in the families ${\mathcal H}_{i(e)}$ and ${\mathcal
  H}_{t(e)}$, where $i(e)$ and $t(e)$ are the initial and
  terminal vertices of $e$.
  \item For any edge $e$, $(G_e,|.|_{{\mathcal H}_e})$ is quasi
  isometrically embedded in $(G_{i(e)},|.|_{{\mathcal
  H}_{i(e)}})$ and in $(G_{t(e)},|.|_{{\mathcal
  H}_{t(e)}})$.
\end{enumerate}
\end{definition}

\begin{remark}
The definition of a graph of strongly relatively hyperbolic groups
is slightly more restrictive than the equivalent definition for
weakly relatively hyperbolic groups. This is because the description
of the subgroups to put in the relative part is heavier in the
former case than in the latter. For the sake of clarity of the
theorem, we adopted Item (b), hoping that this is a not too bad
compromise between clarity and generality.
\end{remark}

We assume given a graph of strongly relatively hyperbolic groups
$\mathcal G$, with fundamental group $\mathcal J$. As before, the
edge and vertex groups are denoted by $G_e$ and $G_v$. Each one
comes with a family of parabolic subgroups, denoted by ${\mathcal
  H}_e$ or ${\mathcal H}_v$.
We construct as before the tree of
hyperbolic spaces $(\widehat{X},{\mathcal T},\pi)$.

\begin{definition}
With the notations above, let $w = v_1 e_1 v_2 \cdots e_k v_{k+1}$
be a $\mathcal T$-geodesic, where the $v_i$'s and $e_i$'s are the
vertices and edges crossed by $w$.

Let $i_{e,v}$ denote the injection of the edge-group
$G_{\overline{e}}$ into the vertex group $G_{\overline{v}}$, where
$\overline{e}$, $\overline{v}$ denote the edge and vertex of
$\mathcal G$ whose lifts contain respectively $e$ and $v$. Then we
denote by $\alpha_w$ the endomorphism with domain a subgroup of
$G_{\overline{v}_1}$ and with image a subgroup of
$G_{\overline{v}_{k+1}}$ given by:
$$\alpha_w = i_{e_k,v_{k+1}} \circ \cdots \circ i_{e_1,v_2} \circ
i^{-1}_{e_1,v_1}$$
\end{definition}

Observe that, if $w$ is a geodesic between two vertices in the lift
of a same vertex of $\mathcal G$, then we can identify $w$ with an
element of the free subgroup of $\mathcal J$ (the fundamental group
of $\mathcal G$) generated by the edges in a complement of a maximal
tree.

\begin{definition}
With the notations above: we say that $H_i \in {\mathcal H}_u$ and
$H_j \in {\mathcal H}_v$ {\em belong to a same orbit of parabolic
subgroups} if there is a geodesic $w$ in $\mathcal T$ from a lift of
$u$ to a lift of $v$ such that $\alpha_w(H_i) = h^{-1} H_j h$ for
some $h \in G_v$.
\end{definition}

There are two kinds of orbits: the finite ones, where all the
$H_i$'s are distinct, and the infinite ones, where infinitely many
conjugates of each $H_i$ appear. An endomorphism $\alpha_w$ between
two subgroups of $G_v$ {\em fixes $H_i$ up to conjugacy} if $H_i$
belongs to the domain of $\alpha_w$ and there is $h \in G_v$ such
that $\alpha_w(H_i) = h^{-1} H_i h$.

\begin{definition}
With the notations above, let $T$ be a maximal tree in $\mathcal G$
and let $F$ be the free subgroup of $\mathcal J$ generated by the
edges in the complement of $T$.

If $H_i \in {\mathcal H}_v$ is a subgroup in an infinite orbit of
parabolic subgroups, then the {\em free extension of $H_i$} is the
subgroup of $\mathcal J$ generated by $H_i$ and by the elements of
the form $a h^{-1}$, where $h \in G_v$ and $a$ is an element in $F$
which satisfies $\alpha_a(H_i) = h^{-1} H_i h$
and which fixes up to conjugacy any parabolic subgroup belonging to an infinite
orbit.
\end{definition}

\begin{remark}
\label{enaijebienbesoin} The subgroup of $F$ which fixes, up to
conjugacy, all the parabolic subgroups in ${\mathcal H}_{v}$ whose
orbit is infinite is finitely generated. For a simple case of free
extension, we refer the reader to Definition \ref{extension}, and to
the particular case of the mapping-torus construction, whose
definition follows \ref{extension}.
\end{remark}

If $g \in G_v$ and $H_i \subset {\mathcal H}_v$, we denote by
$v(gH_i)$ the exceptional vertex of $\widehat{X}$ associated to the
right-class $gH_i$ in the stratum considered.

\begin{definition}
With the notations above, an {\em exceptional leaf} is a maximal set
$S$ of exceptional vertices in $\widehat{X}$ such that: $v(gH_i) \in
S \cap X_{a}$ and
  $v(g^{\prime}H_j)
\in S \cap X_{b}$ if and only if there is $h \in G_{\overline{b}}$
s.t. $\alpha_{[a, b]}(H_i) = h^{-1} H_j h$ and $g^{\prime} =
\alpha_{[a, b]}(g)h^{-1}$.
\end{definition}

\begin{definition}
\label{inutile2} Let $\mathcal G$ be a graph of strongly relatively
hyperbolic groups. The universal covering of $\mathcal G$ satisfies
the {\em strong relative hallways-flare property} if:
\begin{enumerate}
  \item the space $\widehat{X}$ (see above) satisfies the hallways-flare property,
  \item for any $M \geq 0$, there is $T \geq 0$ such that the vertical width of any region
  where two exceptional leaves remain at horizontal distance smaller than $M$
  one from each other is smaller than $T$.
\end{enumerate}
\end{definition}

The second condition in the above definition is needed for the BCP.
Since, by Item (a), the space $\widehat{X}$ satisfies the
hallways-flare property, it suffices in fact that the existence of
$T$ in Item (b) be satisfied for a constant $M$ greater than the
constant of hyperbolicity commonly denoted by this same letter.

 The most general theorem we get is the following one:

\begin{theorem}
\label{je suis content de moi 2} Let $\mathcal G$ be a finite graph
of strongly relatively hyperbolic groups. If the universal covering
of $\mathcal G$ satisfies the strong relative hallways-flare
property, then the fundamental group of $\mathcal G$ is strongly
hyperbolic relative to the family formed by exactly one
representative from each finite orbit of parabolic subgroups and by
the free extensions of exactly one representative from each infinite
orbit of parabolic subgroups.
\end{theorem}

\begin{remark}
For a simple situation of Theorem \ref{je suis content de moi 2},
still giving a good illustration of the phenomena appearing here, we
refer the reader to Theorem \ref{theoreme principal fort}. A simple
example of free extension is given by the mapping-torus of a family
of subgroups, defined after \ref{extension}.
\end{remark}

We begin the proof with the

\begin{lemma}
\label{yoyoyo} With the assumptions and notations of Theorem \ref{je
suis content de moi 2}, there exists $C$ such that any exceptional
leaf is a discrete subset of a $C$-vertical tree.
\end{lemma}

\pols{Since there are only finitely many parabolic subgroups
preserved up to conjugacy and since the free groups which permutes
these subgroups up to conjugacy are finitely generated, there are
only finitely many conjugation elements. Let $m$ be the maximum of
their word-lengths. Then $m+\frac{3}{2}$ ($\frac{1}{2}$ for going
from an exceptional vertex of $\widehat{X}$ to $X$ plus $1$ for
going through a right $H$-class coned in
$\widehat{X}$) gives the announced constant.} \\

The following lemma is a straightforward consequence of the strong
relative hallways-flare property:

\begin{lemma}
\label{blablabla} Two exceptional leaves through two distinct points
in a same stratum of $\widehat{X}$ are connected by a diagonal (see
Definition \ref{def de diagonale}) of horizontal length greater or
equal to $1$, the endpoints of which are exponentially separated in
all the directions outside a region whose vertical size is uniformly
bounded above.
\end{lemma}

\begin{lemma}
\label{gaga} With the assumptions and notations of Theorem \ref{je
suis content de moi 2}, there exists $C(a, b)$ such that, if $g,
g^{\prime}$ are two $(a,b)$-quasi geodesics of ${\widehat{X}}$
between two exceptional leaves  $L_{1}, L_{2}$, then $g, g^{\prime}$
admit decompositions
 $g = g_{1}g_{2}g_{3}$ and $g^{\prime}=g^{\prime}_{1} g^{\prime}_{2}
 g^{\prime}_{3}$ with
the following properties:  $g_{1} \subset \vois{C(a,
b)}{tel}{L_{1}}$, $g^{\prime}_{1} \subset \vois{C(a,
b)}{tel}{L_{1}}$, $g_{3} \subset \vois{C(a, b)}{tel}{L_{2}}$,
$g^{\prime}_{3} \subset \vois{C(a, b)}{tel}{L_{2}}$ and
$d^{H}_{tel}(g_{2}, g^{\prime}_{2}) \leq C(a, b)$. If  $g$ and
$g^{\prime}$ have the same endpoints then $d^{H}_{tel}(g,
g^{\prime}) \leq C(a, b)$.
\end{lemma}

\pols{This is an easy consequence of Theorem \ref{pipeau50}. For
simplicity assume that the attaching-maps of $\widehat{X}$ are quasi
isometries so that Theorem \ref{pipeau} can be applied. The given
two exceptional leaves bound a $C_{\ref{yoyoyo}}$-corridor. Both $g$
and $g^{\prime}$ are approximated by two paths $\mathcal G$ and
${\mathcal G}^\prime$ which only possibly differ by their first and
last maximal vertical segments in $L_1$ and $L_2$. These last
vertical segments are where $g$ and $g^{\prime}$ are not necessarily
close one to each other if they don't have the same endpoints but
are close to the given exceptional leaves. As written before, the
extension to the general case where there is not a corridor, but
only a pseudo-corridor, between the two exceptional leaves, is
easily dealt with by using Theorem \ref{pipeau50} instead of Theorem
\ref{pipeau}.}

\begin{definition}
\label{tjrs+}
We denote by $C(\widehat{X})$ the metric space obtained
from $\widehat{X}$ by putting a cone over each exceptional leaf.
\end{definition}

Lemma \ref{toujoursplus} below stresses the importance of this new
coned space.

\begin{lemma}
\label{toujoursplus} With the assumptions and notations of Theorem
\ref{je suis content de moi 2}: $C(\widehat{X})$ is hyperbolic
and satisfies the BCP with respect to the exceptional leaves if and
only if
the fundamental group of $\mathcal G$ is strongly hyperbolic
relative to the subgroups given by Theorem \ref{je suis content de
moi 2}.
\end{lemma}

\begin{remark}
Assume that $H_1,H_2$ are subgroups of $G$ such that $\alpha(H_1)$
is a conjugate of $H_2$ and $\alpha(H_2)$ is a conjugate of $H_1$.
Then, in $C(\widehat{X})$, cones are put above the right
$H_i$-classes, and their exceptional vertices all belong to a same
exceptional leaf. However, only one of the two subgroups $H_1,H_2$
appears in the subgroups of the relative part described by Theorem
\ref{je suis content de moi 2} because otherwise the condition of
malnormality would be violated.
\end{remark}

\pol{toujoursplus}{Let $\mathcal Y$ be the space obtained by coning
the universal covering of $\mathcal G$ according to the parabolic
subgroups described in Theorem \ref{je suis content de moi 2}. The
essential difference between $\mathcal Y$ and the coned space
$C(\widehat{X})$ of Definition \ref{tjrs+} is the following one:

In $C(\widehat{X})$ a horizontal cone is first put over {\em all} the
right-classes for the parabolic subgroups in the
edge and vertex groups; then a
``vertical cone'' is put over {\em all} the vertices which belong to a
same
exceptional leaf. In $\mathcal Y$, a cone is put on the right-classes
of exactly one
subgroup from each finite orbit, and of exactly one free extension
of subgroup in each infinite orbit.

Observe that in both $C(\widehat{X})$ and  $\mathcal Y$, there is
exactly one exceptional vertex for each exceptional leaf. One thus
has a natural one-to-one correspondence, denoted by $\mathcal B$,
between the exceptional vertices of $C(\widehat{X})$ and those of
$\mathcal Y$.
Assume that there is a horizontal cone in
$C(\widehat{X})$ over two points  $x, y$ in a same stratum of
$\tilde{X}$. It belongs to an exceptional leaf and we denote by
$v(gH)$ the exceptional vertex associated to this leaf. Consider the
exceptional vertex ${\mathcal B}(v(gH))$ of $\mathcal Y$. Assume
that $x, y$ do not belong to the cone with vertex ${\mathcal
B}(v(gH))$. Then there are two points $x^{\prime}, y^{\prime}$ in
another stratum which are at bounded telescopic distance from  $x$
and $y$ and belong to this cone. This is straightforward if $v(gH)$
is the vertex of the cone over a finite exceptional leaf. Otherwise
this comes from the finite generation of the free groups which
preserve the parabolic subgroups up to conjugacy and from the fact
that there is an upper-bound on the length of the conjugacy
elements.

There is a natural map $j \colon \mathcal Y \rightarrow
C(\widehat{X})$ whose restriction to $\tilde{X}$ is the identity-map
and which maps each exceptional vertex $v(gH)$ of $\mathcal Y$ to
the exceptional vertex ${\mathcal B}^{-1}(v(gH))$ of $C(\widehat{X})$.
The observation of the previous paragraph readily implies the
following assertion: if $g$ is a quasi geodesic of $\mathcal Y$,
then $j(g)$ is a quasi geodesic of $C(\widehat{X})$ (with possible
different constants of quasi geodesicity) whose trace in $\tilde{X}$
is Hausdorff-close to the trace of $g$. The lemma follows.}

\begin{remark}
The hyperbolicity of the coned space $C(\widehat{X})$ follows from
the quasi convexity of the exceptional leaves implied by Lemma
\ref{yoyoyo} and from the arguments developed for proving
Proposition 1 of \cite{Sz}. However we re-prove it when listing
below the arguments for checking the BCP.
\end{remark}

\begin{lemma}
\label{ptd} With the notations above: assume that $\widehat{X}$
satisfies the strong hallways-flare property. Let  $g_{1}, g_{2}$ be
two $(a, b)$-quasi geodesics of $C(\widehat{X})$, the terminal
points of which are at most $1$-apart in $\widehat{X}$, and with
same initial point in $\widehat{X}$. Assume the existence of a
generalized corridor $\mathcal C$ between the vertical trees of the
endpoints of $g_{1}$. There exists  $C(a, b, r)$ such that, if the
traces $\widehat{g}_i$'s of the $g_i$'s in $\widehat{X}$ satisfy
$\widehat{g}_{i} \subset \vois{r}{\widehat{X}}{{\mathcal C}}$ for
$i=1,2$ then $d^{H}_{C(\widehat{X})}(g_{1}, g_{2}) \leq C(a, b, r)$.
Furthermore, if $g_{1}$ and $g_{2}$ do not backtrack then they
satisfy the two conditions required for the BCP with a constant
$D(a, b, r)$.
\end{lemma}

We emphasize that this proposition is false if one only requires a
bound on the distance in $C(\widehat{X})$ from the
$g_i$'s to $\mathcal C$. \\

\pols{For simplicity we assume that $\mathcal C$ is a corridor, the
adaptation to generalized corridors is straightforward. We consider
the horizontal projections on $\mathcal C$ of the maximal subpaths
of $g_{1}, g_{2}$ which belong to $\widehat{X}$.
From Lemma
\ref{dimanche}, these projections are $(C_{\ref{dimanche}}(a, b,
r),C_{\ref{dimanche}}(a, b, r))$-quasi geodesics. From Lemmas
\ref{yoyoyo}, \ref{blablabla} on the one hand and Lemma \ref{bo2} on
the other hand, there is $K$, depending on $r$ and $C_{\ref{casi
leafe}}(C_{\ref{yoyoyo}})$, such that the projections of the
exceptional leaves are $K$-vertical trees, for which there exists a
constant $L$ playing the r\^{o}le of the constant
$T_{\ref{inutile2}}$. It is equivalent to prove the announced
properties for the bigon $g_{1},g_{2}$ with respect to the
exceptional leaves than to prove them for the above projections on
$\mathcal C$.

If $g_{1}, g_{2}$ go through the same exceptional leaves, then their
projections on $\mathcal C$ satisfy the same property with respect to
the projections of the exceptional leaves. From Lemma \ref{gaga},
the ``bigon'' obtained by projection to the generalized corridor is
thin. Moreover the points where the projections of  $g_{1}$ and
$g_{2}$ penetrate a given exceptional leaf are close, because either
they are close to the diagonal preceding this exceptional leaf, or
they leave a same exceptional leaf: in this last case we are done by
the existence of the constant  $L$ above (the analog on the corridor
of the constant  $T_{\ref{inutile2}}$). Let us now assume that
$g_{1}$ enters in an exceptional leaf $S$ but $g_{2}$ does not. Of
course this also holds for the respective projections on $\mathcal
C$. We then distinguish three cases:

\noindent {\em First case: the exit point of  $g_{1}$ is
followed by a diagonal with horizontal length greater than some
constant (depending on the constants of hyperbolicity and
exponential separation)}. Then (the projection of) $g_{2}$ has to go
to a bounded neighborhood of this diagonal, this is Theorem
\ref{pipeau}.
It remains before in a bounded
horizontal neighborhood of the exceptional leaf, the bound depending
on $a, b$ and $r$ (since the constants of quasigeodesicity of the
projections depend on  $r$). Thus the vertical length of the passage
of
$g_{1}$ through this exceptional leaf is bounded above by a
constant depending on  $a, b$ and  $r$.

\noindent {\em Second case: the exit point of  $g_{1}$ is
followed by another exceptional leaf}. Thanks to
the existence of the constant $L$ and Lemma \ref{blablabla}, we can
follow the same arguments as above, appealing to Proposition
\ref{proposition importante 1} rather than directly Theorem \ref{pipeau}. We
leave the reader
work out details and computations.

\noindent {\em Third case: the exit point of $g_{1}$ is followed by a
  horizontal geodesic with horizontal length bounded above by the
  constant
of the first case.} In this case, this horizontal geodesic ends
at the vertical boundary  of $\mathcal C$. The entrance-point of
$g_{1}$ in  $S$ is close to a point in  $g_{2}$. Since $g_{2}$ is a  $(a, b)$-quasi geodesic and
$g_{2}$ does
not pass through  $S$, it cannot happen that the passage of $g_{1}$
though $S$ is a long passage at small horizontal distance from the
considered vertical boundary. Thus, if it is a long passage, then
there is a stratum,  which is nearest to the entrance-point of
$g_{1}$ in  $S$, where the horizontal distance between  $S$ and the
considered vertical boundary is smaller than the critical
constant. From Proposition \ref{proposition importante 1},  $g_{2}$
lies in a bounded neighborhood of $S$ until reaching this
stratum. Once again, this gives an upper-bound on the vertical length
of $S$.

The proof of Lemma \ref{ptd} now follows in an easy way: to
conclude for the BCP, we need of course the fact that the horizontal
metrics on the strata satisfy the BCP.}

\begin{proposition}
\label{un peu complique} With the assumptions of Lemma
\ref{ptd}: there exist $C(a,b) \geq 1$ and $D(a,b)
> 0$ such that, if $x_0,x_1,\cdots,$ $x_n$ are consecutive points
in some exceptional leaf $L$, which lie outside the horizontal
$D(a,b)$-neighborhood of a generalized corridor $\mathcal C$, and if
the vertical distance between the strata of $x_0$ and $x_n$ is
greater than $C(a,b)$, then no non-backtracking $(a,b)$-quasi
geodesic of $C(\widehat{X})$ with both endpoints in the horizontal
$D(a,b)$-neighborhood of $\mathcal C$ contains as subpath the cone
over $\{x_{0},x_{n}\}$.
\end{proposition}

See proof in subsection \ref{upc} of Section \ref{quasi convexe}.\\

\pot{je suis content de moi 2}{Let $g,g^\prime$ be two
non-backtracking $(a,b)$-quasi geodesics of $C(\widehat{X})$ with
same initial point, and with terminal points at most $1$-apart in
$\widehat{X}$. We assume for a while that the attaching-maps of
$\widehat{X}$ are quasi isometries. There is thus a corridor
$\mathcal C$ between vertical trees passing through the initial and
terminal points of $g$.
\par Let $p$ be a passage of  $g$ (resp. of $g^{\prime}$) through the
cone over a subset $S$ of an exceptional leaf
outside
the $D_{\ref{un peu complique}}(a,b)$-neighborhood of $\mathcal C$
in $\widehat{X}$. From Proposition \ref{un peu complique},
substituting  $p$  by  $S$ yields
a non-backtracking $(\kappa(a, b),\kappa^{\prime}(a,b))$-quasi
geodesics  $h$ (resp. $h^{\prime}$) of $C(\widehat{X})$, with $\kappa(a, b) =
C_{\ref{un peu complique}}(a,b) C_{\ref{yoyoyo}}a$ and
$\kappa^{\prime}(a, b) = C_{\ref{un peu complique}}(a, b)
C_{\ref{yoyoyo}} (b+1)$, such that $d^{H}_{C(\widehat{X})}(g, h) \leq
1$
(resp.  $d^{H}_{C(\widehat{X})}(g^{\prime},
h^{\prime}) \leq 1$). We can thus assume that all passages like $p$ have
been suppressed in $h$ and $h^{\prime}$ as above.

By Proposition \ref{quasiconvexite}, the subpaths of $h$ and
$h^{\prime}$ between two exceptional
 leaves are contained in the
 horizontal $C_{\ref{quasiconvexite}}(\kappa(a, b),$ $\kappa^{\prime}(a,
 b))$-neighborhood of a corridor between these leaves. Thus $h$ and $h^\prime$
 are contained in the $D_{\ref{un peu complique}}(a,b) +
 C_{\ref{quasiconvexite}}(\kappa(a, b),\kappa^{\prime}(a,
 b))$-neighborhood of ${\mathcal C}$ in $\widehat{X}$. From Lemma
 \ref{ptd}, $h,h^\prime$ satisfy the BCP. The conclusion for
 $g,g^\prime$ follows.

 The proof of the hyperbolicity follows the same scheme. If
 $g,g^\prime$ form a $(a,b)$-quasi geodesic bigon of
 $C(\widehat{X})$, one first substitutes it by a
 non-backtracking $(a,b)$-quasi geodesic bigon $g_0,g^\prime_0$
 with $d^H_{C(\widehat{X})}(g,g_0) \leq b$,
 $d^H_{C(\widehat{X})}(g^\prime,g^\prime_0) \leq b$. The line of the
 arguments thereafter is the same than above: at the end, Lemma
 \ref{ptd}
gives the thinness
 of the quasi geodesic bigons instead of the BCP. As in Section \ref{hts},
 the hyperbolicity follows
 from Lemma \ref{ilya gautero}.

 Let us briefly sketch the adaptations of the above arguments to
 the general case, where the attaching-maps are not necessarily quasi isometries.
 Instead of a corridor, there is a pseudo-corridor between two given
 points (see Definition \ref{pseudo corridor}).
 This is a concatenation of generalized
corridors ${\mathcal C}_0,\cdots,{\mathcal C}_n$ connected by
horizontal geodesics $h_1,\cdots,h_n$. From this decomposition, $g$
and $g^\prime$ inherit a decomposition $g_0 \cdots g_n$, $g^\prime_0
\cdots g^\prime_n$ such that the initial and terminal points of
$g_i, g^\prime_i$ lie in the $(2 \delta + 1)$-neighborhood, in
$\widehat{X}$, of ${\mathcal C}_i$. Up to increasing the constant
$b$ to a constant $b^\prime = b + 2 \delta + 1$, and extending a
little bit $g_i$ and $g^\prime_i$, we can assume that each pair
$g_i,g^\prime_i$ forms a $(a,b^\prime)$-quasi geodesic bigon.
Proving the BCP (resp. the thinness) for any of the
$g_i,g^\prime_i$, is enough to prove the BCP (resp. the thinness)
for $g,g^\prime$. We so go back to the beginning of the proof.}

 \subsection{Some corollaries of Theorem \ref{je suis content de moi 2}}

 \label{scott}

 \begin{definition}
\label{extension} Let $G$ be a finitely generated group and let
${\mathcal H} = \{H_{1}, \cdots, H_{k}\}$ be a finite family of
subgroups of $G$. Let $\F{r} = \langle {\mathcal A} \rangle$, with
${\mathcal A} = \{\alpha^{\pm 1}_1,\cdots,\alpha^{\pm 1}_r\}$, be a
rank $r$ free group of relative automorphisms of $(G,{\mathcal H})$.
A {\em $\F{r}$-extension
 ${\mathcal H}_{\mathcal A}$} of $\mathcal H$ is a maximal family
of subgroups of $G \rtimes_{\mathcal A} \F{r}$ of the form $\langle
H_{i}, a_{i,1} g^{-1}_{i,1}, \cdots,a_{i,m} g^{-1}_{i,m},\cdots
\rangle$ such that:

\begin{itemize}
  \item each $a_{i,m} \in \F{r}$ satisfies $a_{i,m} H_i = g^{-1}_{i,m} H_i g_{i,m}$
and $ \langle a_{i,1},a_{i,2},\cdots \rangle$ generates the subgroup
of all the elements of $\F{r}$ which fixes each $H_i$ up to
conjugacy;
  \item if ${\mathcal H}_j, {\mathcal H}_{j^\prime}$ are two distinct
subgroups in ${\mathcal H}_{\mathcal A}$ with $H_i \in {\mathcal
H}_j$, $H_k \in {\mathcal H}_{j^\prime}$ then no element of
$\langle {\mathcal A} \rangle$ conjugates $H_i$ to $H_k$ in $G
\rtimes_{\mathcal A} \F{r}$.
\end{itemize}
\end{definition}

When $r=1$ in the above definition, i.e. $\F{r} = \langle t
\rangle$, we get the easier notion of the {\em mapping-torus of
$\mathcal H$} under a relative
  automorphism $\alpha$ of $(G,{\mathcal H})$. This
  is a maximal family ${\mathcal H}_\alpha$ of subgroups
  ${\mathcal H}_j \subset G_\alpha$ satisfying the following properties:
  \begin{itemize}
    \item ${\mathcal H}_j = \langle H_{i_j},t^{n_{i_j}} g^{-1}_{i_j} \rangle$,
  where $n_{i_j}$ is the minimal integer such that there is $g_{i_j} \in G$
  with $\alpha^{n_{i_j}}(H_{i_j}) = g^{-1}_{i_j} H_{i_j} g_{i_j}$;
    \item whenever ${\mathcal H}_j = \langle H_{i_j},t^{n_{i_j}} g^{-1}_{i_j} \rangle$,
    ${\mathcal H}_k = \langle H_{i_k},t^{n_{i_k}} g^{-1}_{i_k} \rangle$
  are two distinct subgroups in ${\mathcal H}_\alpha$, no power of $t$ conjugates $H_{i_j}$ to
  $H_{i_k}$ in $G_\alpha$.
  \end{itemize}

\begin{remark}
As was previously noticed in Remark \ref{enaijebienbesoin}, the free
subgroup of $\F{r}$ which fixes each $H_i$ up to conjugacy is
finitely generated (see \cite{GH} for a proof).
\end{remark}

 \begin{theorem}
\label{theoreme principal fort} Let $G$ be a finitely generated
group and let $\mathcal H$ be
 a finite family of infinite subgroups  of $G$. Let $\F{r}$ be a
uniform free group of relatively hyperbolic automorphisms of $(G,
{\mathcal H})$. Then, if $G$ is strongly hyperbolic
relative to $\mathcal H$, $G \rtimes \F{r}$
 is strongly hyperbolic relative to the $\F{r}$-extension of $\mathcal H$.
\end{theorem}

\pots{It suffices to check that the definition of a uniform free
group of relatively hyperbolic automorphisms implies the {\em
strong} hallways-flare property. The exponential separation of the
vertical segments is clear but one has to prove that {\em any two}
exceptional leaves also separate exponentially one from each other.
Assume that this is not satisfied. Then, there is $M \geq 0$ such
that for any $N \geq 1$, there is $\alpha_{w} \in \F{r}$ with $|w|
\geq N$, s.t. there is a geodesic word $u$ in $(G,|.|_{{\mathcal
H}})$ of the form $h_{1} H_{i_{1}} h_{2} \cdots H_{i_{k}}
 h_{k+1}$ (where  $h_{j}$ stands for a passage of the geodesic
in the Cayley graph of $G$ whereas $H_{i_{j}}$ stands for a passage
of the geodesic in a right-class for $H_{i_j}$) satisfying the
following properties:

\begin{enumerate}
  \item $|u|_{{\mathcal H}} \leq M$,
 \item the image under $\alpha_w$ of any element with geodesic word $HuH^\prime$ has the form
 $rHuH^\prime s$, where $H,H^\prime$ stand for passages through right-classes
 for the corresponding parabolic subgroups, and where the relative lengths of $r$ and $s$ only depend on the
 length of $w$.
\end{enumerate}

Here $H$ and $H^\prime$ are the parabolic subgroups of $G$
corresponding to the right-classes associated to the two exceptional
leaves which violate, for the considered $w$, the strong exponential
separation property. The existence of $u$ above comes from the
finiteness of the family $\mathcal H$ and from the finite generation
of $G$: they imply together that there are only finitely many
geodesic words of a given form which have relative length smaller
than $M$.

Since $G$ is strongly hyperbolic relative to $\mathcal H$, $\mathcal
H$ is almost malnormal in $G$. This readily implies, by choosing
elements in $H$ and $H^\prime$ which are sufficiently long enough in
$(G,|.|_{S})$, that there is an element $g$ of the form $HuH^\prime
. H^\prime u^{-1} H \equiv HuH^\prime u^{-1}H$ which is not
conjugate to an element of a parabolic subgroup. Furthermore $g$ can
be chosen not to be a torsion element. From Corollary 4.20 of
\cite{Os}, $\displaystyle \lim_{n \to + \infty} |g^n|_{{\mathcal H}}
= + \infty$. However $\alpha_w(g)$ has the form $rHuH^\prime s
s^{-1} H^\prime u^{-1} H r^{-1} \equiv r HuH^\prime u^{-1}H r^{-1}$.
Thus $|\alpha_w(g^n)|_{{\mathcal H}} \leq |g^n|_{{\mathcal H}} + 2
|r|_{{\mathcal H}}$. Since $|r|_{{\mathcal H}}$ is a constant only
depending on $|w|_{{\mathcal H}}$, by choosing $n$ sufficiently
large enough we get a contradiction with the uniform hyperbolicity
of $\F{r}$.}

\begin{definition}
A finite family $\{H_1,\cdots,H_k\}$ of subgroups of a group $G$ is
{\em almost malnormal} if:
\begin{enumerate}
  \item for any $i=1,\cdots,k$, $H_i$ is almost malnormal in $G$.
  \item for any $i,j \in \{1,\cdots,k\}$ with $i \neq j$,
the cardinality of the set $\{w \in H_{j} \mbox{ ; } \exists g \in G
\mbox{ s.t. } w \in g^{-1} H_i g\}$ is finite.
\end{enumerate}
\end{definition}

If the family of subgroups consists of only one subgroup, the
definition above is nothing else than the definition of almost
malnormality of this subgroup. It is now widely known that a hyperbolic group $G$ is strongly
hyperbolic relative to any almost malnormal finite family of quasi convex subgroups.
As a corollary of the previous
theorem we thus have:

\begin{corollary}
\label{premier corollaire fort} Let $G$ be a hyperbolic group, let
$\mathcal H$ be a finite family of infinite subgroups of $G$ and let
$\alpha$ be an automorphism of $G$ which is hyperbolic relative to
$\mathcal H$. If $\mathcal H$ is quasi convex and malnormal in $G$
then the mapping-torus group $G_{\alpha} = G \rtimes_\alpha \mz$ is
strongly hyperbolic relative to the mapping-torus of $\mathcal H$.
\end{corollary}

This corollary may be specialized to torsion free one-ended
hyperbolic groups, and so in particular to fundamental groups of
surfaces. We so re-prove the result of \cite{Ga2}. Since there we
gave only an idea for the statement and the proof in the Gromov
relative hyperbolicity case, we include here the full statement of
this result:

\begin{corollary}
\label{rappel0} Let $G$ be a torsion free one-ended hyperbolic group
and let $\alpha$ be an automorphism of $G$. Let $\mathcal H$ be a
maximal family of maximal subgroups of $G$ which consist entirely of
elements on which $\alpha$ acts up to conjugacy periodically or with
linear growth. Then $G_{\alpha}$ is weakly hyperbolic relative to
$\mathcal H$, and strongly hyperbolic relative to the mapping-torus
of $\mathcal H$.

If $G$ is the fundamental group of a compact surface  $S$ (possibly
with boundary) with negative Euler characteristic and $h$ a
homeomorphism of $S$ inducing  $\alpha$ on  $\pi_{1}(S)$ (up to
inner automorphism), then the subgroups in $\mathcal H$ are:
\begin{enumerate}[(i)]
  \item the cyclic subgroups generated by the boundary curves,
  \item the subgroups associated to the maximal subsurfaces which are unions
of components on which $h$ acts periodically, pasted together along
reduction curves of the Nielsen-Thurston decomposition,
  \item the cyclic
subgroups generated by the reduction curves not contained in the
previous subsurfaces.
\end{enumerate}
\end{corollary}

\pocs{From Corollary \ref{premier corollaire fort}, we only have to
prove that the considered automorphism $\alpha$ of $G$ is hyperbolic
relative to the given family of subgroups. The passage from the
surface case to the torsion free one-ended hyperbolic group case is
done thanks to the JSJ-decomposition theorems of \cite{Bo1}. We
refer the reader to \cite{Ga2} for more precisions and concentrate
on the surface case. The fundamental group of $S$ is the fundamental
group of a graph of groups $\mathcal G$ such that:
\begin{itemize}
  \item the edge groups are cyclic
subgroups associated to the reduction curves and boundary
components,
  \item the vertex groups are the subgroups associated to the pseudo-Anosov components (type $I$ vertices) and
  to the maximal subsurfaces with no pseudo-Anosov components (type $II$ vertices),
  \item the (outer) automorphism $\alpha$ induced by the homeomorphism preserves the graph of
  groups structure.
\end{itemize}
We consider the universal covering of $\mathcal G$ and the
associated tree of spaces. We measure the length of a geodesic in
this tree of spaces as follows:

\begin{itemize}
  \item we count zero for the passages through the edge-spaces and
through the type $II$ vertex-spaces,
  \item we measure the length of the pieces through the type $I$ vertex-spaces by integrating against
  the stable and unstable measures of the invariant foliations (a boundary-component is considered to
  belong to both invariant foliations and so the contribution of a path in such a leaf amounts to zero).
\end{itemize}

There is $N \geq 1$ such that, when the total stable (resp.
unstable) length of a geodesic in a type $I$-vertex space is two
times its unstable (resp. stable) length, then it is dilated by a
factor $\lambda > 1$ under $N$ iterations of $\alpha^{-1}$ (resp. of
$\alpha$). In the other cases, we find $N \geq 1$ such that the
total length is dilated under $N$ iterations of both $\alpha$ and
$\alpha^{-1}$. Similar computations have been presented in
\cite{Ga2}. The conclusion of the relative hyperbolicity of $\alpha$
now comes easily since pieces with positive length, dilated either
under $\alpha^N$ or under $\alpha^{-N}$, and pieces with zero length
alternate.} \\

 Up to now, we only exhibited extensions of relatively hyperbolic
groups via semi-direct products. However such a product is only a
particular case of HNN-extension. Alibegovic in \cite{Ali}, Dahmani
in \cite{Da} or Osin in \cite{Os2} treat acylindrical HNN-extensions
and amalgated products. Let us now give a theorem about
non-acylindrical HNN-extensions. Theorem \ref{injectif} below deals
with injective, not necessarily surjective, endomorphisms of
relatively hyperbolic groups. We first introduce a notion of
relative malnormality.

\begin{definition}
Let $G$ be a group and let ${\mathcal H} = \{H_1,\cdots,H_k\}$ be a
finite family of subgroups of $G$. A subgroup $H^\prime$ of $G$ is
{\em almost malnormal relative to $\mathcal H$} if there is an
upper-bound on the  $\mathcal H$-word length of the elements in the
set $\{w \in H^{\prime}  \mbox{ ; } \exists g \in G - H^{\prime}
\mbox{ with } w \in g^{-1} H^\prime g\}$.
\end{definition}

If $\mathcal H$ is empty, the definition above is nothing else than
the usual notion of almost malnormality and if in addition there is
no torsion, we get the notion of malnormality. Whereas the
definitions of a relative automorphism and of a mapping-torus of a
family of subgroups given in Definition \ref{premieres definitions}
remain valid for injective endomorphisms, the definition of relative
hyperbolicity for automorphisms is easily adapted to the more
general case of injective endomorphisms:

\begin{definition}
Let $G$ be a finitely generated group and let $\mathcal H$ be a
finite family of subgroups of $G$. An injective endomorphism
$\alpha$ of $G$ is {\em hyperbolic relative to $\mathcal H$} if
$\alpha$ is a relative endomorphism of $(G, {\mathcal H})$ and there
exist $\lambda > 1$ and $M,N \geq 1$ such that, for any $w \in
\mbox{Im}(\alpha^{N})$ with $|w|_{{\mathcal H}} \geq M$, if
$|\alpha^N(w)|_{\mathcal H} \geq \lambda |w|_{\mathcal H}$ does not
hold then $w = \alpha^N(w^\prime)$ with $|w^\prime|_{\mathcal H}
\geq \lambda |w|_{\mathcal H}$.
\end{definition}

\begin{theorem}
\label{injectif} Let $G$ be a finitely generated group, let $\alpha$
be an injective endomorphism of $G$ and let $G_{\alpha}$ be the
associated mapping-torus group, i.e. the associated ascending
HNN-extension. Let $\mathcal H$ be a finite family of infinite
subgroups of $G$ such that $\alpha$ is hyperbolic relative to
$\mathcal H$. Assume that $\mathrm{Im}(\alpha)$ is almost malnormal
relative to $\mathcal H$. Then, if $G$ is strongly hyperbolic
relative to $\mathcal H$, $G_{\alpha}$ is weakly hyperbolic relative
to $\mathcal H$ and strongly hyperbolic relative to the
mapping-torus of $\mathcal H$.
\end{theorem}

\begin{remark}
The reader will notice at once that the above theorem does not treat
the extension of weakly relatively hyperbolic groups. The reason is
that the condition of relative almost malnormality does not imply in
this case the relative hallways-flare property. This last property
is however also a necessary condition, although we do not give here
a direct proof: in the absolute hyperbolicity case, Gersten was the
first to give the converse to the combination theorem, using
cohomological arguments \cite{Ge} and we adapt his arguments in
\cite{GH}. Bowditch exposed a more direct proof in \cite{Bo2}.
\end{remark}

\pot{injectif}{We first prove the following

 \begin{lemma}
\label{aieaieaieaie} Let $G$ be a finitely generated group which is
strongly hyperbolic relative to a finite family of subgroups
$\mathcal H$. Let $K$ be a finitely generated subgroup of $G$, which
is almost malnormal relative to $\mathcal H$, which is strongly
hyperbolic relative to a (possibly empty) finite family ${\mathcal
H}^\prime$ the subgroups of which are conjugated to subgroups in
${\mathcal H}$, and such that $(K, |.|_{{\mathcal H}^\prime})$ is
quasi isometrically embedded in $(G,|.|_{\mathcal H})$. There exists
$C > 0$ such that, if $x, y$ (resp.  $z, t$) are any two vertices in
a same right-class $gK$ (resp.  $h K$) with $g \neq h$ then
$d_{{\mathcal H}}(P_{[z, t]}(x), P_{[z, t]}(y)) \leq C$.
\end{lemma}

\pols{Since $(G,|.|_{\mathcal H})$ is hyperbolic, there is a
constant $\delta \geq 0$ such that the geodesic triangles of
$(G,|.|_{\mathcal H})$ are $\delta$-thin. Thus, geodesic rectangles
are $2 \delta$-thin. This implies the existence of a quadruple of
vertices $x_0,y_0,z_0,t_0$ with $x_0,y_0 \in [x,y]$, $z_0,t_0 \in
[z,t]$ and $d_{\mathcal H}(x_0,z_0) \leq 2 \delta + 1$, $d_{\mathcal
H}(y_0,t_0) \leq 2 \delta +1$. Since $(K,|.|_{{\mathcal H}^\prime})$
is $(\lambda,\mu)$-quasi isometrically embedded in
$(G,|.|_{{\mathcal H}})$, and $(G,|.|_{{\mathcal H}})$ is
$\delta$-hyperbolic, there exist $c_0(\lambda,\mu,\delta)$ and
$x_1,y_1,z_1,t_1$ such that $g^{-1}x_1, g^{-1} y_1 \in K$, $h^{-1}
z_1, h^{-1} t_1 \in K$ and $d_{{\mathcal H}}(x_0,x_1) \leq
c_0(\lambda,\mu,\delta)$, $d_{{\mathcal H}}(y_0,y_1) \leq
c_0(\lambda,\mu,\delta)$, $d_{{\mathcal H}}(z_0,z_1) \leq
c_0(\lambda,\mu,\delta)$, $d_{{\mathcal H}}(t_0,t_1) \leq
c_0(\lambda,\mu,\delta)$. We choose $x_1,y_1,z_1,t_1$ to minimize
the distance {\em in $(G,S)$} (that is the distance associated to
the given finite set of generators $S$ of $G$) respectively to
$x_{0}, y_{0}, z_{0}, t_{0}$. We denote by $[x_1,y_1]_K$ (resp.
$[z_1,t_1]_K$) the images, under the embedding of $K$ in $G$, of
geodesics between the pre-images of $x_1,y_1$ (resp. $z_1,t_1$) in
$K$. Both $[x_1,y_1]_K$ and $[z_1,t_1]_K$ are $(\lambda,\mu)$-quasi
geodesics. Moreover $[x_1, z_1] [z_1,t_1]_K [t_1,y_1]$ is a
$(\lambda,4 \delta + 2 + 4 c_0(\lambda,\mu,\delta) + \mu)$-quasi
geodesic between $x_1$ and $y_1$. Since $G$ is strongly hyperbolic
relative to $\mathcal H$, $G$ satisfies the BCP with respect to
$\mathcal H$. This gives a constant $c_1(\lambda,\mu,\delta)$ such
that the $\mathcal H$-classes $[x_1,z_1]$ and $[t_1,y_1]$ go through
correspond to geodesics in $(G,S)$ with length smaller than
$c_1(\lambda,\mu,\delta)$: indeed, since  $x_{1}, y_{1}, z_{1},
t_{1}$ were chosen to minimize the distances in $(G,S)$ with respect
to $x_{0}, y_{0}, z_{0}, t_{0}$, the $\mathcal H$-classes crossed by
$[x_1,z_1]$ and $[t_1,y_1]$ are not crossed by $[x_{1}, y_{1}]_{K}$.
Therefore the distance in $(G,S)$ between $x_1$ and $z_1$ on the one
hand, and between $y_1$ and $t_1$ on the other hand is less or equal
to $(2 \delta + 1 + 2 c_0(\lambda,\mu,\delta))
c_1(\lambda,\mu,\delta)$. There are a finite number of elements in
$G$ with such an upper-bound on the length, measured with a
word-metric associated to a finite set of generators. Whence, by the
almost normality of $K$ relative to $\mathcal H$, an upper-bound on
the length between $x_1$ and $y_1$, and so also between $x_0$ and
$y_0$. Lemma \ref{aieaieaieaie} is proved.} \\ \par From Lemma
\ref{aieaieaieaie}, the overlapping of two distinct right
$\mathrm{Im}(\alpha)$-classes is bounded above by a constant.
Together with the fact that $\alpha$ is a relatively hyperbolic
endomorphism, this implies the exponential separation property.
Getting the strong version of this property is done as in the proof
of Theorem \ref{theoreme principal fort}. Theorem \ref{injectif} now
follows from Theorem \ref{je suis content de moi 2}.}

\section{Proof of Proposition \ref{proposition importante 1}}

\label{propimpo1}

\noindent {\small \bf Convention:} {\em Throughout the paper, the
constants of hyperbolicity and of quasi isometry are chosen
sufficiently large enough to satisfy the conclusions of Lemma
\ref{indispensable}, and also sufficiently large enough so that
computations make sense. Moreover the horizontal subpaths of the
$(a,b)$-quasi geodesics considered will be assumed to be horizontal
geodesics.} The hyperbolicity of the strata gives a constant
$C(a,b)$ such that any $(a,b)$-quasi geodesic  $g$ may be
substituted by another one $g^{\prime}$ with $d^{H}_{tel}(g,
g^{\prime}) \leq C(a,b)$
and satisfying this latter property. \\

Our first lemma is about quasi geodesics. It holds not only in a
corridor but in the whole tree of hyperbolic spaces.

\begin{lemma}
\label{deuxieme debut} Let $(\tilde{X},{\mathcal T},\pi)$ be a tree
of hyperbolic spaces with exponentially separated $v$-vertical
segments. Let $g$ be a $v$-telescopic $(a,b)$-quasi geodesic in
$\tilde{X}$. There exist $C(a,b)$ and $D$ such that, if $[x,y]
\subset g \cap X_{\alpha}$ satisfies $d_{hor}(x,y) \geq C(a,b)$ then
for any ${\mathcal T}$-geodesic $w$ starting at $\alpha$ with
$|w|_{\mathcal T} \geq D + n t_0$, $n \geq 1$, we have
$d^{i}_{hor}(wx,wy) \geq \lambda^{n} d_{hor}(x,y)$.
\end{lemma}

\pols{We denote by  $\lambda > 1, M, t_{0} \geq 1$ the constants
  of hyperbolicity and by $\lambda_{+}, \mu$ the constants of quasi isometry.
  Let us choose $n_\star(a)$ such that $\frac{a}{\lambda^{n_\star}}
< 1$. Solving the inequality $e
> a (\frac{1}{\lambda^{n_\star}} e + 2 n_\star t_0) + b$ gives us
$e(a,b) \geq \frac{2 a n_\star t_0 + b}{1 - a
\frac{1}{\lambda^{n_\star}}}$.

{\em \noindent Claim:} If $d_{hor}(x,y) \geq e(a,b)$, if
$x^\prime,y^\prime$ are the endpoints of two  $v$-vertical segments
$s,s^\prime$ of vertical length $n_\star t_0$, starting at $x$ and
$y$ and with  $\pi(s)=\pi(s^{\prime})$, then for any
${\mathcal T}$-geodesic $w_{0}$
such that $w_0\pi(s)$ is a ${\mathcal T}$-geodesic and
$|w_0|_{\mathcal T} = t_0$,
$d^{i}_{hor}(w_0x^\prime,w_0y^\prime)
\geq \lambda d_{hor}(x^\prime,y^\prime)$ holds.

{\em \noindent Proof of Claim:} Assume the existence of $w$ with
$|w|_{\mathcal T} = n_\star t_0$ such that for some
$x^\prime,y^\prime$ with $x \in w x^\prime$, $y \in w y^\prime$ and
$d_{hor}(x^\prime,y^\prime) \geq M$, $d_{hor}(x,y) \geq
{\lambda^{n_\star}} d_{hor}(x^\prime,y^\prime)$ holds. Then
$\frac{1}{\lambda^{n_\star}} e + 2 n_\star t_0$ is the telescopic
length of a telescopic path between $x$ and $y$. But the inequality
given at the beginning of the proof tells us that the existence of
such a telescopic path is a contradiction with the fact that $g$ is
a  $v$-telescopic $(a,b)$-quasi geodesic. Therefore, if
$d_{hor}(x,y) \geq e(a,b)$ and $d_{hor}(x,y) \geq
\lambda^{n_\star}_+ (M+\mu)$ (this last inequality is to assert that
$d_{hor}(x^\prime,y^\prime) \geq M$ - see above), then
$d_{hor}(x^\prime,y^\prime)$ does not increase after $t_0$ in the
direction of the  $v$-vertical segments $s,s^\prime$. The claim
follows from the exponential separation of the $v$-vertical
segments.
\par From the inequality given by the Claim, since
$d_{hor}(x^\prime,y^\prime) \geq \lambda^{-n_\star}_+
(d_{hor}(x,y)+\mu)$, we easily compute an integer $N_\star$ such
that, if $w_0$ is as in the Claim but with length $N_\star t_0$ then
$d^{i}_{hor}([w_0\pi(s)]x,[w_0\pi(s)]y) \geq \lambda d_{hor}(x,y)$.
Setting $D = N_\star t_0$ and $C(a,b) = e(a,b)$, the constant
computed above, we get the lemma.}
\\

{\noindent \small \bf Notations:} $\delta$ a fixed non negative
constant, $(\tilde{X},{\mathcal T},\pi)$ a tree of
$\delta$-hyperbolic spaces, $\mathcal C$ a generalized corridor with
exponentially separated $v$-vertical segments, $\lambda > 1 ,M,t_0
\geq 1$ the associated constants of hyperbolicity, $\lambda_+, \mu$
the associated constants of quasi isometry, $g$ a $v$-telescopic
$(a,b)$-quasi geodesic of $\mathcal C$. The above constants are
chosen sufficiently large enough to satisfy the conclusions of Lemma
\ref{casi leafe}.

\begin{lemma}
\label{the end peutetre}  There exists $C(a,b)$ such that, if the
endpoints  $x,y$ of $g$ both lie in a same stratum $X_\alpha$, if $d_{hor}(x,y) \geq
C(a,b)$ then, for any ${\mathcal T}$-geodesic $w$ starting at
$\alpha$ with $|w|_{\mathcal T} \geq C(a, b) + n t_0$, $n \geq 1$,
and $w \cap \pi(g) = \{\alpha\}$,  we have:

$$d^{i}_{hor}(wx,wy) \geq \lambda^n d_{hor}(x,y).$$
\end{lemma}

\pols{Let us observe that, if $[p, q]$ is any
  horizontal geodesic in $g$ then the $v$-vertical trees of
  $p$ and $q$ bound a horizontal geodesic  $[p^{\prime}, q^{\prime}]$ in
  $[x,y]$.

\noindent{\em Claim:} If $d_{hor}(p^\prime,q^\prime) \geq
  Cte$, with $Cte \equiv \lambda^{t_0}_{+}
  (C_{\ref{deuxieme debut}}(a,b)+t_0+\mu)$, then for any $w$ as given by the current
      Lemma with
   $|w|_{\mathcal T} \geq D_{\ref{deuxieme debut}} + t_0$,
$d^{i}_{hor}(wp^{\prime},wq^{\prime}) \geq \lambda
d_{hor}(p^\prime,q^\prime)$.

\noindent{\em Proof of Claim:} If $p^\prime$ and $q^\prime$ are not
exponentially separated in the direction of $p,q$ after $t_0$, then,
because of the hallways-flare property, they are exponentially
separated after $t_0$ in the direction of $w$, which yields the
announced inequality. Let us assume that $p^\prime,q^\prime$ are
separated after $t_0$ in the direction of $[\pi(p^\prime),\pi(p)]$.
Thus $d^i_{hor}(rp^\prime,rq^\prime) \geq \lambda^n
d_{hor}(p^\prime,q^\prime)$ for ${\mathcal T}$-geodesic $r$ with
$|r|_{\mathcal T} = nt_0$ and $r \cap w = \{\alpha\}$. Therefore
$d_{hor}(p,q) \geq C_{\ref{deuxieme debut}}(a,b)+t_0$. Lemma
\ref{deuxieme debut} then implies that $p,q$ are exponentially
separated in the direction of $[\pi(p),\pi(p^\prime)]$ after
$D_{\ref{deuxieme debut}}+t_0$, and the claim is proved.

There is a finite decomposition of $[x,y] \subset X_\alpha$ in
subgeodesics $[p^\prime_j,q^\prime_j]$ with disjoint interiors such
that each $[p^\prime_j,q^\prime_j]$ connects two $v$-vertical trees
through the endpoints of a maximal horizontal geodesic in $g$. We
denote by $I_{D}$ the set of $[p^\prime_j,q^\prime_j]$'s with
$d_{hor}(p^\prime_j,q^\prime_j) \geq Cte$ and by $I_{C}$ the set of
the others. Let us choose an integer $n \geq 1$. We consider a
stratum $X_\beta$ with $d_{\mathcal T}(\beta,\alpha) =
D_{\ref{deuxieme debut}} + n t_{0}$. Let $h$ be the horizontal
geodesic in ${\mathcal C} \cap X_\beta$ which connects the two
$v$-vertical trees through $x$ and $y$. Assume that the endpoints of
$h$ are exponentially separated after $t_{0}$ in the direction of
$[\beta,\alpha]$. Then:

\begin{equation}
\label{mon equation} \lambda^{n} |I_{D}|_{hor} \leq |h|_{hor} \leq
\lambda^{-n} (|I_{D}|_{hor} + |I_{C}|_{hor})
\end{equation}

so that $$|I_{C}|_{hor} \geq
\frac{\lambda^{n}-\lambda^{-n}}{\lambda^{-n}} |I_{D}|_{hor}$$ and
consequently, since $d_{hor}(x,y) = |I_D|_{hor} + |I_C|_{hor}$,
$$|I_{C}|_{hor} \geq \frac{X(n)}{1+X(n)} d_{hor}(x,y)$$ with $X(n)
= \frac{\lambda^{n}-\lambda^{-n}}{\lambda^{-n}}$. Since
$\displaystyle \lim_{n \to + \infty} \frac{X(n)}{1+X(n)} = 1$,
there is $n_{\star} \geq 0$ such that for any $n \geq n_{\star}$,
$$|I_{C}|_{hor} \geq \frac{1}{2} d_{hor}(x,y).$$ But, by
definition, the horizontal length of each subgeodesic in $I_{C}$ is
smaller than $Cte$. Thus the number of elements in $I_C$ is at least
the integer part of $\frac{1}{2 Cte} d_{hor}(x,y) +1$. Furthermore,
since $g$ is a $v$-telescopic path, the telescopic length of any
subpath of $g$ containing $j$ maximal horizontal geodesics is at
least $(j-1)$. We so obtain:
$$|g|^{v}_{tel} \geq \frac{1}{2 Cte} d_{hor}(x,y).$$
On the other hand: $$d^{v}_{tel}(x,y) \leq \lambda^{-n}
d_{hor}(x,y) + 2 n t_0.$$ since there is a $v$-telescopic path
between $x$ and $y$ the telescopic length of which is given by the
right-hand side of the above inequality. Since $g$ is a
$(a,b)$-quasi geodesic, the last two inequalities give $n_{\star
\star} \geq 0$ such that for $n \geq n_{\star \star}$:
$$d_{hor}(x,y) \leq \frac{2ant_0+b}{\frac{1}{2 Cte} -
a \lambda^{-n}}.$$ Taking the maximum of $n_\star,n_{\star \star}$
and the above upper-bound for $d_{hor}(x,y)$, we get the announced
constant in the case where the endpoints of the horizontal geodesic
$h$ above are exponentially separated in the direction of
$[\beta,\alpha]$. If not, there are in all the other directions so
that we easily get a constant $N \geq 0$ such that
$d^{i}_{hor}(wx,wy) \geq \lambda d_{hor}(x,y)$ for any ${\mathcal
  T}$-geodesic $w$ with $|w|_{\mathcal T} = N t_0$ and
$[\pi(x),\pi(h)] \subset w$. Lemma \ref{the end peutetre} is then
easily deduced.}
\\

As a consequence we have:

\begin{corollary}
\label{troisieme lemme} With the assumptions and notations of Lemma
\ref{the end peutetre}, there exists $C(a,b,d) \geq d$ such that, if
$x,y$ are the endpoints of two $v$-vertical segments $s,s^\prime$
with $d^{i}_{hor}(s,s^\prime) \leq d$ with  $\pi(s) = \pi(s^{\prime})$
and
 $\pi(s) \cap \pi(g) = \{\alpha\}$, then $d_{hor}(x,y) \leq C(a,b,d)$.
\end{corollary}

\begin{remark}
At this point, we would like to notice that Lemma \ref{the end
peutetre} is similar to Lemma 6.7 of \cite{Ga}. However in addition
of some misprints, a slight mistake took place there in the proof of
the Lemma. Indeed the inequality (\ref{mon equation}) in the proof
of Lemma \ref{the end peutetre} is true here, in the generalized
corridor, but there the constant $\lambda$ should have been modified
to take into account the so-called ``cancellations''.
\end{remark}

\begin{lemma}
\label{quatrieme lemme} Let $x$ and $y$ be the endpoints of a
$r$-vertical segment $s$ in  $\mathcal C$. There exists $C(r)$ such
that, if the intersection-point $z$ of a $v$-vertical tree through
$y$ with the stratum $X_{\pi(x)}$ satisfies $d_{hor}(x,z) \geq
C(r)$, then for any ${\mathcal T}$-geodesic $w$ with $|w|_{\mathcal
T} = nt_0$,  $n \geq 1$, and $w \cap \pi(s) = \{\pi(x)\}$, $d^{i}_{hor}(wx,wz) \geq
\lambda^{n} d_{hor}(x,z)$.
\end{lemma}

\pols{If $|s|_{vert} \leq t_0$, the existence of the constants of
quasi isometry, Item (a) of Lemma \ref{indispensable2}, and the
definition of a $r$-vertical segment give an upper-bound for
$d_{hor}(x,z)$. Let us thus assume $|s|_{vert} > t_0$. Choose $d$
such that $\lambda d - r^\prime \geq 2r^\prime$, where $r^\prime$ is
the above upper-bound when $|s|_{vert} = t_0$. Then set $C =
\max(d,M)$. Assume $d_{hor}(x,z) \geq C$ and that $x$ and $z$ are
exponentially separated in the direction given by $s$. If
$[\pi(x),\pi(y)] = w_0 w^\prime$ with $|w_0|_{\mathcal T} = t_0$,
then $d^{i}_{hor}(w_0x,w_0z) \geq \lambda d_{hor}(x,z)$. Thanks to the
inequality used to defined  $d$, one easily concludes
that the horizontal distance between $s$ and the vertical tree
through $y$ increases along $s$ when going from $x$ to $y$ which of
course cannot happen. The conclusion follows from the hallways-flare
property.} \\

\pop{proposition importante 1}{We set $x_i = T_i \cap X_\alpha$, in
particular $d_{hor}(x_1,x_2)=L$. We consider the region $R$ with
vertical width $C_{\ref{troisieme lemme}}(a,b,L)$ centered at the
stratum $X_\alpha$. We decompose $\mathcal G$ in three subpaths: the
first one, denoted ${\mathcal G}_0$, from the initial point of
$\mathcal G$ until the first point $z$ in ${\mathcal G} \cap R$, the
second one, denoted ${\mathcal G}_1$, from $z$ to the last point $t$
in ${\mathcal G} \cap R$, the third one, denoted ${\mathcal G}_2$,
from $t$ to the terminal point of $\mathcal G$. Obviously the
subpath ${\mathcal G}_1$ can be approximated by the concatenation of
two vertical segments with a horizontal geodesic in $X_\alpha$ (the
approximation constant only depend on $L,a$ and $b$). We denote by
${\mathcal G}^\prime_1$ the resulting path.

We now consider a maximal subpath in ${\mathcal G}_0$ which
satisfies the following properties:

\begin{itemize}
  \item its endpoints lie in a same stratum $X_\beta$,
  \item its image under $\pi$ does not intersect $[\alpha,\beta]$,
  outside $\beta$.
\end{itemize} \par From Corollary \ref{troisieme lemme}, the endpoints of such a
subpath are at horizontal distance smaller than $C_{\ref{troisieme
lemme}}(a,b,L)$ one to each other. Thus, by substituting each such
subpath by a horizontal geodesic connecting its endpoints, we
construct a $C_{\ref{troisieme lemme}}(a,b,L)$-vertical segment
${\mathcal G}^\prime_0$. We do the same for ${\mathcal G}_2$, so
obtaining a $C_{\ref{troisieme lemme}}(a,b,L)$-vertical segment
${\mathcal G}^\prime_2$. From Lemma \ref{quatrieme lemme},
${\mathcal G}^\prime = {\mathcal G}^\prime_0 \cup {\mathcal
G}^\prime_1 \cup {\mathcal G}^\prime_2$ lies in a bounded
neighborhood of the $v$-vertical segments connecting its endpoints
to $x_1$ and $x_2$. From the construction, $d^H_{tel}({\mathcal
G},{\mathcal G}^\prime) \leq a C_{\ref{troisieme lemme}}(a,b,L) + b
+ 1$. The proposition follows.}

\section{Quasiconvexity of corridors}
\label{quasi convexe}

In this section we prove Propositions \ref{quasiconvexite}, \ref{qc}
and \ref{un peu complique}.

\subsection{Two basic lemmas}

We need first a very general lemma about Gromov hyperbolic spaces.

\begin{lemma}
\label{reecriture} Let $(X,d)$ be a Gromov hyperbolic space. There
exists an increasing
affine function $D(r) \geq 0$, and $C \geq 0$
such that, if $[x,y]$ is a diameter of a ball $B_{x_0}(r)$ with $r
\geq C$ and $w$ is any path in $X$ with $w \cap B_{x_0}(r) =
\{x,y\}$, then $|w|_d \geq e^{D(r)}$.
\end{lemma}

\noindent This lemma is a rewriting of Lemma 1.6 page 26 of
\cite{CDP}. \hfill $\Box$

\begin{lemma}
\label{quadrilatere} Let $\tilde{X}$ be a tree of
$\delta$-hyperbolic spaces which satisfies the hallways-flare
property. There exists $C$ such that, if $x,y,z,t$ are the vertices
of a geodesic quadrilateral in some stratum $X_\alpha$, with
$d_{hor}(x,z) \leq 2 \delta$, $d_{hor}(y,t) \leq 2 \delta$, and
$d_{hor}(x,y) \geq C$, $d_{hor}(z,t) \geq C$, then for any ${\mathcal
  T}$-geodesic $w$ with $|w|_{\mathcal T} \geq
C_{\ref{indispensable}} + nt_0$, starting at $\pi(x)$:

$$d^{i}_{hor}(wx,wy) \geq \lambda^n d_{hor}(x,y) \Leftrightarrow
d^{i}_{hor}(wz,wt) \geq \lambda^n d_{hor}(z,t)$$
\end{lemma}

\pols{If $A,B$ are two subsets of a metric space $(X,d)$, we set
$\displaystyle d^{s}(A, B) = \sup_{x \in A, y \in B} d(x, y)$. Let
us consider any ${\mathcal T}$-geodesic $w$ with
 $|w|_{\mathcal T} = t_{0}$
starting at $\alpha$. From Lemma \ref{indispensable2},
 $$d^{s}_{hor}(wx, wz) \leq \lambda^{t_{0}}_{+} (2 \delta+\mu)$$  and
 $$d^{s}_{hor}(wy, wt) \leq \lambda^{t_{0}}_{+} (2 \delta+\mu).$$ Assume
$d^{i}_{hor}(wx,wy) \geq \lambda d_{hor}(x,y)$ but
$d^{i}_{hor}(wz,wt) < \lambda d_{hor}(z,t)$.

We take $d_{hor}(x, y) \geq M$
 and  $d_{hor}(z, t) \geq M$. Assume $d^{s}_{hor}(wz, wt) \leq \frac{1}{\lambda} d_{hor}(z,
 t)$. But $d_{hor}(z, t) \leq 4 \delta + d_{hor}(x, y)$. Putting
 together these inequalities we get
 $$\lambda d_{hor}(x,y) \leq 2  \lambda^{t_{0}}_{+} (2 \delta+\mu) +
\frac{1}{\lambda} (4 \delta + d_{hor}(x, y)).$$ Whence an upper
bound for $d_{hor}(x,y)$ and thus for  $d_{hor}(z, t)$. If
$d^{s}_{hor}(wz, wt) > \frac{1}{\lambda} d_{hor}(z,
 t)$ then the lemma follows from the definition of the constant
 $C_{\ref{indispensable}}$, see the corresponding lemma.} \\

The above two lemmas are not needed if one only considers trees of
$0$-hyperbolic spaces, the proof in this last case being much
simpler.

\subsection{Approximation of quasi geodesics with bounded vertical
  deviation} \hfill

\label{aoqgwbvd}

Lemma \ref{oubli} below states that, in a tree of hyperbolic
spaces $(\tilde{X}, {\mathcal T})$, a quasi geodesic with bounded
image in $\mathcal T$ lies close to a corridor between its
endpoints. This is intuitively obvious and nothing is new neither
surprising in the arguments of the proof: they heavily rely upon the
$\delta$-hyperbolicity of the strata and the fact that strata are
quasi isometrically embedded into each other. For the sake of
brevity, we do not develop them here.

\begin{lemma}
\label{oubli} Let  $(\tilde{X}, {\mathcal T}, \pi)$ be a tree of
hyperbolic spaces. There exists $C(\kappa, a, b)$ such that, if $g$
is any $v$-telescopic $(a,b)$-quasi geodesic with
$\mathrm{diam}_{{\mathcal T}}(\pi(g)) \leq \kappa$, if $\mathcal C$
is a generalized corridor between its endpoints then $g \subset
\vois{C(\kappa, a, b)}{tel}{{\mathcal C}}$.
\end{lemma}

\subsection{Stairs}

\label{st}

The sign $\simeq_1$ stands for an equality up to $\pm 1$.

\begin{definition}
\label{escaliers fibre} Let $\mathcal C$ be a generalized corridor
in a tree of hyperbolic spaces $(\tilde{X},{\mathcal T},\pi)$.

A {\em $r$-stair relative to $\mathcal C$}, $r \geq M$, is a
telescopic path $\mathcal S$ the maximal vertical segments of which
have
  vertical length greater than $C_{\ref{indispensable}}$ and such that, for any maximal horizontal
geodesic $[a_{i}, b_{i}]$ in $\mathcal S$:
\begin{enumerate}
  \item  $d_{hor}(a_{i}, b_{i}) \geq r$ and $d^{i}_{hor}([a_{i}, b_{i}], {\mathcal C})
  \simeq_1
    d_{hor}(a_{i},P^{hor}_{{\mathcal C}}(a_i))$,
  \item any pair of points $a,b \in [a_i,b_i]$ with
  $d_{hor}(a,b) \geq r$ are exponentially separated in
  the direction of
  the $\mathcal T$-geodesic $[\pi(a_{i}), \pi(a_{i+1})]$.
\end{enumerate}
\end{definition}

\begin{lemma}
\label{stairs fibre} With the notations of Definition \ref{escaliers
fibre}: there exist  $C \geq C_{\ref{quadrilatere}}$ such that for
any $r \geq C$, if $\mathcal S$ is a $r$-stair relative to $\mathcal
C$, if $\mathcal U$ is a generalized corridor between a vertical
tree through the terminal point of $\mathcal S$ and a vertical
boundary of $\mathcal C$, then
$${\mathcal S} \subset \vois{r+2\delta}{hor}{{\mathcal U}}.$$
\end{lemma}

\pols{Let $a_{i}, b_i \in {\mathcal S}$ as given in Definition
\ref{escaliers fibre} and let $z$ be a point at the intersection of
the stratum $X_{\pi(a_i)}$ with a vertical tree through some point
farther in the stair. Then:

\noindent {\em Claim 1:} There exists $K > 0$ not depending on $a_i$
nor $z$ such that, if $r$ is sufficiently large enough then
$d^{i}_{hor}([a_{i}, z], {\mathcal C}) \geq
d_{hor}(a_{i},P^{hor}_{{\mathcal C}}(a_i)) - K$.

\noindent {\em Proof of Claim 1:} Choose $K$ such that
$e^{D_{\ref{reecriture}}(K)} > 4 \delta + 1$ and assume
$d^{i}_{hor}([a_{i}, z], {\mathcal C}) <
d_{hor}(a_{i},P^{hor}_{{\mathcal C}}(a_i)) - K$. Then Lemma
\ref{reecriture} implies that $[b_{i}, z]$ descends at least until a
$2 \delta$-neighborhood of $a_{i}$. Assume $r \geq
C_{\ref{quadrilatere}}+2 \delta$. Then Lemma \ref{quadrilatere}
gives an initial segment of $[b_{i},z]$ of horizontal length greater
than $r - 2 \delta$ which is dilated in the direction of
$[\pi(a_{i}), \pi(a_{i+1})]$. If $r$ is chosen sufficiently large
enough with respect to the constants of hyperbolicity for a corridor
(see Lemma \ref{bo2}), we get $z^{\prime}$ at the intersection of
the considered vertical tree through $z$ with the stratum
$X_{\pi(a_{i+1})}$ such that $d^{i}_{hor}([a_{i+1}, z^{\prime}],
{\mathcal C}) < d_{hor}(a_{i+1},P^{hor}_{{\mathcal C}}(a_{i+1})) -
K$. The repetition of these arguments show that the horizontal
distance between $\mathcal S$ and the vertical tree through $z$ does
not decrease along $\mathcal S$. This is an absurdity since $z$ was
chosen in a vertical tree through a point farther in $\mathcal S$.
The proof of Claim 1 is complete.

\noindent {\em Claim 2:} There exists $K(r)$ not depending on $b_i$
nor $z$ such that, if $r$ is sufficiently large enough then
$d^{i}_{hor}([b_{i}, z], {\mathcal C}) \geq
d_{hor}(b_i,P^{hor}_{{\mathcal C}}(b_i)) - K(r)$.

\noindent {\em Proof of Claim 2:} Let $z_\star \in [b_i,z]$ with
$d_{hor}(z_\star,P^{hor}_{{\mathcal C}}(z_\star)) \simeq_1
\max(d^i_{hor}([b_i,z],{\mathcal C}),d_{hor}(a_i,$
$P^{hor}_{{\mathcal C}}(a_i)))$. From the $\delta$-hyperbolicity of
the strata, $[b_i,z_\star]$ lies in the horizontal $2
\delta$-neighborhood of $[a_i,b_i]$. Assume $d_{hor}(b_i,z_\star)
\geq r$ and is sufficiently large enough to apply Lemma
\ref{quadrilatere}. Then there is $K(r)$ such that, if $z_\star$
satisfies $d_{hor}(z_\star,P^{hor}_{{\mathcal C}}(z_\star)) <
d_{hor}(b_i,P^{hor}_{{\mathcal C}}(b_i)) - K(r)$, the points $b_i$
and $z_\star$ are exponentially separated in the direction of
$[\pi(a_i),\pi(a_{i+1})]$. We thus obtain at $a_{i+1}$ a situation
similar to that of Claim 1. The proof of Claim 2 follows.

Lemma \ref{stairs fibre} is easily deduced from the above two
claims, we leave the reader work out the easy details.}

\begin{lemma}
\label{esperons} There exists $C > 0$ such that, for any $r \geq
C_{\ref{stairs fibre}}$, if $\mathcal S$ is a $r$-stair relative to
$\mathcal C$, which is not contained in the vertical
$C$-neighborhood of the stratum containing its initial point, then
the terminal point of $\mathcal S$ does not belong to the telescopic
$r$-neighborhood of $\mathcal C$.
\end{lemma}

\pols{Decompose $\mathcal S$ in maximal
  substairs  ${\mathcal S}_{0} \cdots {\mathcal S}_{k}$ such that
  $\pi({\mathcal S}_{j})$ is a geodesic of $\mathcal T$.
  Let $[a_i,b_i]$ be the first maximal horizontal geodesic in ${\mathcal
S}_j$, let $x$ be the initial point of ${\mathcal S}_j$ and let $z$
be any point in ${\mathcal S}_j$ with $n t_0 \leq d_{\mathcal
T}(\pi(z),\pi(x)) \leq (n+1) t_0$.

The inequality
\begin{equation}
\label{pupu} d_{hor}(z,P^{hor}_{{\mathcal C}}(z)) \geq Cte
\lambda^{n} d_{hor}(a_i,b_i)
\end{equation}
is an easy consequence of the definition of a stair and of Lemma
\ref{quadrilatere} as soon as $r \geq C_{\ref{quadrilatere}}$.
Indeed, the initial segment of horizontal length $r$ in
$[b_i,P^{hor}_{\mathcal C}(b_i)]$ lies in the horizontal $2
\delta$-neighborhood of $[b_i,a_i]$. The assertion then follows from
Item (b) of Definition \ref{escaliers fibre} and Lemma
\ref{quadrilatere}.

The inequality (\ref{pupu}) readily gives the announced result.}

\subsection{Approximation of a quasi geodesic by a stair} \hfill

\label{popneufun}

\noindent {\small \bf Notations:}  $(\tilde{X}, {\mathcal T})$ a
tree of  $\delta$-hyperbolic spaces with exponentially separated
$v$-vertical segments, $v \geq C_{\ref{preparation}}$,  $\mathcal C$
a generalized corridor, $g$ a $v$-telescopic $(a,b)$-quasi geodesic.

\begin{lemma}
\label{technique} Assume that the endpoints of $g$ are in a
horizontal $r$-neighborhood of $\mathcal C$ and that $g$ lies in the
closed complement of this horizontal neighborhood. Suppose moreover
that the maximal vertical segments in $g$ have vertical length
greater than $3 (C_{\ref{indispensable}} + D_{\ref{deuxieme
debut}})$.

Then there exist $C(r,a,b), D(a,b), E(r,a,b)$ such that for any $r
\geq D(a,b)$, either $g$ lies in the telescopic
$C(r,a,b)$-neighborhood of a $E(r,a,b)$-stair relative to $\mathcal
C$, where $E(r,a,b)$ is affine in $r$, or $g$ is contained in the
telescopic $C(r, a, b)$-neighborhood of $\mathcal C$.
\end{lemma}

\pols{We decompose the proof in two steps. The first one is only a
warm-up, to present the ideas in a particular, but important, case.
The general case, detailed in the second step, is technically more
involved
but no new phenomenon appears. \\

\noindent{\em Step $1$: Proof of Lemma \ref{technique} when the
horizontal length of any maximal horizontal subpath in $g$ is
greater than some constant (depending on $a$ et $b$).} The endpoints
of any horizontal subpath $h$ of $g$ with horizontal length greater
than $C_{\ref{deuxieme debut}}(a,b)$ are exponentially separated
under every geodesic $w$ of $\mathcal T$ with length
$D_{\ref{deuxieme debut}}$. If $|h|_{hor} \geq
C_{\ref{quadrilatere}}$, this is also true for any horizontal
geodesic $h^\prime$ in the $2 \delta$-neighborhood of $h$. Finally,
if $|h|_{hor}$ is sufficiently large enough, by Lemma \ref{bo2} the
endpoints of $h$ are also exponentially separated in any
$v$-corridor containing $h$. If $e(a,b)$ is the maximum of the above
constants, we now assume $|h|_{hor} \geq 3 e(a,b)$.

Let us consider two consecutive maximal horizontal geodesics
$h_1,h_2$ in $g$, separated by a vertical segment $s$. Let $\mathcal
D$ be a corridor containing $h_1$ and $s$. Then:

\begin{equation}
\label{jen rajoute} |h_2 \cap \vois{2 \delta}{hor}{{\mathcal
D}}|_{hor} \leq e(a,b).
\end{equation}

Otherwise we have a contradiction with the fact that the endpoints
of any subgeodesic of $h_2$ whose length is greater than
$C_{\ref{deuxieme debut}}(a,b)$ are exponentially separated in the
direction of $h_1$. \par From the inequality (\ref{jen rajoute}),
the concatenation of $h_1$, $s$ and $h_2$ is $e(a,b)$-close, with
respect to the horizontal distance, of a $2 e(a,b)$-stair relative
to $\mathcal C$ if $d^i_{hor}(h_1,{\mathcal C}) \simeq_1
d_{hor}(a_1,P^{hor}_{{\mathcal C}}(a_1))$ where $a_1$ is the initial
point of $h_1$.

Let us now set $r \geq 3  e(a,b)$ and assume that the maximal
horizontal geodesics in $g$ have horizontal length greater than $r$.
Let $x$ be the initial point of $g$ (in particular
$d_{hor}(x,P^{hor}_{{\mathcal C}}(x)) \simeq_1 r$). Let $s$ be the
vertical segment starting at $x$ and ending at $y$ in $g$. Let $h$
be the maximal horizontal geodesic following $s$ along $g$. Let $n
\geq 1$ be the greatest integer with $n (C_{\ref{indispensable}} +
D_{\ref{deuxieme debut}}) \leq |s|_{vert}$.

By assumption $x$ and $P^{hor}_{{\mathcal C}}(x)$ are exponentially
separated in the direction of $s$. Since the strata are quasi
isometrically embedded one into each other, this gives $\kappa > 1$
such that, any two points $a,b \in [x,P^{hor}_{{\mathcal C}}(x)]$
with $d_{hor}(a,b) \geq \mathrm{max}(\frac{1}{\kappa} r,M)$ satisfy
$d_{hor}(\pi(s)a,\pi(s)b) \geq \lambda^n d_{hor}(a,b)$. Thus the
same arguments as those exposed above when working with $h_1, h_2$
show that $|h \cap \vois{2 \delta}{hor}{[y,P^{hor}_{{\mathcal
C}}(y)]}|_{hor} \leq \mathrm{max}(e(a,b),\frac{1}{\lambda^n \kappa}
r,M)$. If $n$ is greater than some critical constant $n_*$, this
last maximum is equal to $e(a,b)$. Thus, in this case, we can take
$h_1 = [x,P^{hor}_{{\mathcal C}}(x)]$ and $h_2 = h$: the above
arguments prove that the concatenation of $h_1$, $s$ and $h_2$ is
$e(a,b)$-close to a $e(a,b)$-stair. If $n$ is smaller than $n_*$,
then we substitute $r$ by $\lambda^{n_* (C_{\ref{indispensable}} +
D_{\ref{deuxieme debut}})}_+ r$, modify $g$ by taking the starting
point at the endpoint $y$ of $s$ and take $h_1$ as the first maximal
horizontal geodesic.

In both cases, by repeating the arguments above at any two
consecutive maximal horizontal geodesic following the first two ones
along $g$, we show that $g$ is $e(a,b)$-close, with respect to the
horizontal distance, of a
$e(a,b)$-stair relative to $\mathcal C$. \hfill $\Box$ \\

\noindent{\em Step $2$: Adaptation of the argument to the general
case:} The boundary trees of $\mathcal C$ are denoted by $L_1$ and
$L_2$, and $g$ goes from $L_1$ to $L_2$. We choose a positive
constant $r$, which when necessary will be set sufficiently large
enough with respect to the constants $C_{\ref{stairs fibre}}$, $M,
\delta$ and $C_{\ref{quadrilatere}}$. Let $x_{0}$ be the initial
point of $g$. It lies in the boundary of the horizontal
$r$-neighborhood of $\mathcal C$. We denote by ${\mathcal C}_{i}$
and $x_{i}$, $i=1,\cdots$, a sequence of corridors and points of $g$
defined inductively as follows:
\begin{enumerate}
  \item ${\mathcal C}_{i}$ is a corridor with boundary trees a $v$-vertical
  tree through $x_{i-1}$ and the $v$-vertical boundary $L_2$ of ${\mathcal C}$,
  \item $x_{i}$ is the first point following $x_{i-1}$ along  $g$ such
that $d_{hor}(x_{i},P^{hor}_{{\mathcal C}_{i}}(x_i)) \geq r$.
\end{enumerate}
The subpath of $g$ between $x_{i-1}$ and $x_{i}$ is denoted by
$g_{i-1, i}$. Obviously $g_{i-1,i}$ is contained in the horizontal
 $r$-neighborhood of ${\mathcal C}_{i}$. We project it to ${\mathcal
 C}_{i}$. From Lemma \ref{dimanche}, we get a $C_{\ref{casi
   leafe}}(v)$-telescopic $(C_{\ref{dimanche}}(a,b,r),
C_{\ref{dimanche}}(a,b,r))$-quasi geodesic of  $({\mathcal
  C}_{i},d^{C_{\ref{casi leafe}}(v)}_{tel})$. We set $X(a,b,r) =
C_{\ref{proposition importante 1}}(r, C_{\ref{dimanche}}(a,b,r),
C_{\ref{dimanche}}(a,b,r))$. From Proposition \ref{proposition
importante 1}, $P^{hor}_{{\mathcal C}_{i}}(g_{i-1,i})$ is contained
in the $X(a,b,r)$-neighborhood of the concatenation of a subpath of
$[x_{i-1}, P^{hor}_{{\mathcal
    C}_{i-1}}(x_{i-1})]$ with a vertical segment in  ${\mathcal
  C}_{i}$ (and is followed by $[P^{hor}_{{\mathcal C}_{i}}(x_{i}),
x_{i}]$). Consider in this approximation of (a subpath of) $g$ a
maximal collection of points $y_{i}$ which defines a $r$-stair
relative to ${\mathcal C}$. The points $y_{i}$ do not necessarily
agree with the $x_{i}$'s, because it might happen that, after
$x_{i-1}$ for instance, the approximation constructed above reenters
in the $r$-neighborhood of  ${\mathcal C}_{{i-1}}$ before leaving
the $r$-neighborhood of ${\mathcal C}_{i}$. We proceed as in Step
$1$ and choose the $y_i$'s so that:

\begin{enumerate}
  \item either $y_i$ is contained in a maximal horizontal geodesic,
  and from the observations in Step $1$, this horizontal geodesic
  may be included in a stair,
  \item or the vertical distance from $y_i$ to the next horizontal
  geodesic is at least $C_{\ref{indispensable}} + D_{\ref{deuxieme
  debut}}$.
\end{enumerate}

Either we obtain a non-trivial $r$-stair relative to ${\mathcal C}$
which approximates a subpath $g^{\prime}_{0}$ of $g$ or the
approximation we constructed above exhausts $g$ and is contained in
some telescopic neighborhood of ${\mathcal C}$ the size of which is
obtained from the previously exhibited constants. In this last case,
the same assertion holds for the whole path $g$. This is one of the
announced alternatives.

We can thus assume that we got $y_{0}, \cdots, y_{k}$ forming a
$r$-stair relative to $\mathcal C$. It is denoted by $S$. Since the
strata are quasi isometrically embedded one into each other, there
is $\kappa > 1$, only depending on the constants of quasi isometry,
such that $S$ is in fact a $\mathrm{max}(\frac{1}{\kappa}
r,M,e(a,b))$-stair relative to $\mathcal C$. As soon as $r > \kappa
(M+e(a,b))$, which we suppose from now, this maximum is just
$\frac{1}{\kappa} r$. Thus $S$ is a $\frac{r}{\kappa}$-stair whose
maximal horizontal geodesics have horizontal length at least $r$.

By construction $S$ approximates $g^{\prime}_{0} \subset g$. We now
consider the maximal subpath $g^{\prime}_{1}$ of $g$ starting at (or
near - recall that we constructed an approximation of a subpath of
$g$) $y_{k}$ which lies in the $r$-neighborhood of ${\mathcal
C}_{k}$. This last corridor plays the r\^{o}le of the corridor
${\mathcal U}$ of Lemma \ref{stairs fibre}. We project the subpath
$g^{\prime}_{1}$ to ${\mathcal C}_k$, so getting a
$(C_{\ref{dimanche}}(a,b,r), C_{\ref{dimanche}}(a,b,r))$-quasi
geodesic of this corridor. From Lemma \ref{stairs fibre}, and
because of the hyperbolicity of the strata, each horizontal geodesic
of the $\frac{r}{\kappa}$-stair $S$ admits a subgeodesic with
horizontal length greater than  $\frac{\kappa-1}{\kappa} r$ in the
horizontal $2 \delta$-neighborhood of ${\mathcal C}_{k}$. If $r$ is
chosen sufficiently large enough, Lemma \ref{quadrilatere} gives
horizontal geodesics in ${\mathcal C}_{k}$ with horizontal length
greater than $M$ which are dilated in the same directions than the
horizontal geodesics of $S$. Now Proposition \ref{proposition
importante 1} applies and allows us to approximate the projection of
$g^{\prime}_{1}$ on ${\mathcal C}_{k}$ by a sequence of these
horizontal geodesics. But each one of these horizontal geodesics is
close to a point in $g^{\prime}_{0} \subset g$. Thus, since $g$ is a
$(a, b)$-quasi geodesic, the vertical length of $g^{\prime}_{1}$,
and so its telescopic length, is bounded above by a constant
depending on $a$ and $b$. So we can forget $g^{\prime}_{1}$ and
continue the construction of our $\frac{r}{\kappa}$-stair relative
to $\mathcal C$ at the point where the approximation of
$g^{\prime}_{1}$ leaves the $r$-neighborhood of ${\mathcal C}_{k}$.
We eventually exhaust $g$ and obtain a $\frac{r}{\kappa}$-stair
relative to $\mathcal C$. }

\subsection{Proof of Proposition \ref{quasiconvexite}}

Let $g$ and $\mathcal C$ be as given by this proposition. Assume
that some subpath $g^{\prime}$ of $g$ leaves and then reenters the
horizontal $D_{\ref{technique}}(a,b)$-neighborhood of $\mathcal C$.
Assume that $g^\prime$ is not contained in the telescopic
$C_{\ref{technique}}(D_{\ref{technique}}(a,b),a,b)$-neighborhood of
$\mathcal C$.

Suppose for the moment that the vertical segments in $g^\prime$ have
vertical length greater than $3 (C_{\ref{indispensable}} +
D_{\ref{deuxieme debut}})$. Then Lemma \ref{technique} gives
$\mathcal G$, a $R(a,b)$-stair relative to $\mathcal C$, where
$R(a,b) \equiv E_{\ref{technique}}(D_{\ref{technique}}(a,b),a,b)$,
with $d^H_{tel}(g^\prime,{\mathcal G}) \leq
C_{\ref{technique}}(D_{\ref{technique}}(a,b),a,b)$. From Lemma
\ref{esperons}, $\mathcal G$ does not leave the vertical
$C_{\ref{esperons}}(R(a,b))$-neighborhood of the stratum containing
the initial point of $\mathcal G$. Therefore, setting  $$V(a,b) =
C_{\ref{esperons}}(R(a,b)) +
C_{\ref{technique}}(D_{\ref{technique}}(a,b),a,b),$$ $g^\prime$ does
not leave the vertical $V(a,b)$-neighborhood of this stratum. From
Lemma \ref{oubli}, $g^\prime$ lies in the telescopic
$C_{\ref{oubli}}(V(a,b),a,b)$-neighborhood of $\mathcal C$.

It remains to consider the case where the vertical segments in
$g^\prime$ are not sufficiently large enough. Let $s$ be a vertical
segment in $g$ with $|s|_{vert} < 3 (C_{\ref{indispensable}} +
D_{\ref{deuxieme debut}}) \equiv X$.

$({\dag})$ Thanks to the assumption that all the attaching-maps of
the tree of hyperbolic spaces are quasi isometries, $s$ is contained
in a vertical segment $s^\prime$ of vertical length greater than
$X$. We modify $g^\prime$ by sliding, along $s^\prime$, a horizontal
geodesic in $g^\prime$ incident to $s$ until getting a vertical
segment with vertical length $X$. This yields a new telescopic
$(a^\prime,b^\prime)$-quasi geodesic in a bounded neighborhood of
$g$, where the constants $a^\prime,b^\prime$ only depend on $a,b$
and on the constants of quasi isometry. After finitely many such
moves, we obtain a quasi geodesic as desired, and we are done. Since
the vertical distance between two strata is uniformly bounded away
from zero, after finitely many such substitutions, we eventually get
a quasi geodesic, in a bounded neighborhood of $g$, which satisfies
the assumptions required by Lemma \ref{technique}. This completes
the proof of Proposition \ref{quasiconvexite}. \hfill $\Box$

\subsection{Proof of Proposition \ref{qc}}

\label{popqc}

We leave to the reader the usual modifications to pass from
corridors to generalized corridors. There remains the problem of
getting a telescopic path with maximal vertical segments
sufficiently large enough. We start from the sentence marked by a
$(\dag)$ in the preceding subsection. If $s$ is not contained in a
vertical segment $s^\prime$ of vertical length greater than $X$, we
obtain a vertical segment ${\bf s}$ from $b_i$ to $a_{i+1}$
satisfying the following properties (we still denote by $g^\prime$
the $(a^\prime,b^\prime)$-quasi geodesic eventually obtained, we
denote by ${\bf s}_0$ the vertical segment of $g^\prime$ ending at
$a_i$ and by ${\bf s}_1$ the one starting at $b_{i+1}$):

\begin{enumerate}
  \item there is no vertical segment starting at $a_i$ (resp. at $a_{i+1}$) over the
  edge $\pi({\bf s})$ (resp. over $\pi({\bf s}_1)$);
  \item there is no vertical segment ending at $b_i$ over $\pi({\bf s}_0)$.
\end{enumerate}

Consider horizontal geodesics $\alpha_i = [a_i,P^{hor}_{{\mathcal
C}}(a_i)]$, $\beta_i = [b_i,P^{hor}_{{\mathcal C}}(b_i)]$,
$\alpha_{i+1} = [a_{i+1},$ $P^{hor}_{{\mathcal C}}(a_{i+1})]$ and
$\beta_{i+1} = [b_{i+1},P^{hor}_{{\mathcal C}}(b_{i+1})]$. By the
$\delta$-hyperbolicity of the strata, there is $a^\prime_i \in
[a_i,b_i] \cap \vois{2 \delta}{hor}{\alpha_i \cup \beta_i}$ and
$b^\prime_i \in [a_{i+1},b_{i+1}] \cap \vois{2
\delta}{hor}{\alpha_{i+1} \cup \beta_{i+1}}$. Because the strata are
quasi isometrically embedded one into each other, we get two points
$a^{\prime \prime}_i, b^{\prime \prime}_i$ which satisfy:

\begin{enumerate}
  \item[{\rm (A)}] they are $Y$-close (with respect to the horizontal distance)
respectively to $a^\prime_i$ and $b^\prime_i$, where the constant
$Y$ only depends on $\delta$ and on the constants of quasi isometry;
  \item[{\rm (B)}] there is a $v$-vertical segment from $a^{\prime \prime}_i$
  to $b^{\prime \prime}_i$ which is contained in a larger $v$-vertical
  segment going over $\pi({\bf s}_0)$ and $\pi({\bf s}_1)$.
\end{enumerate}

We modify $g^\prime$ by going from $a_i$ to $a^{\prime \prime}_i$
then to $b^{\prime \prime}_i$ and eventually end at $b_{i+1}$. The
resulting path is a $(a^{\prime \prime},b^{\prime \prime})$-quasi
geodesic, where the constants $a^{\prime \prime}, b^{\prime \prime}$
only depends on $\delta$ and on the constants of quasi isometry.
Moreover this new path is in a bounded neighborhood of $g^\prime$.
Thanks to Item (B), we can modify it by enlarging the vertical
segment from $a^{\prime \prime}_i$ to $b^{\prime \prime}_i$. The
conclusion in then the same as in the preceding subsection. \hfill
$\Box$

\subsection{Proof of Proposition \ref{un peu complique}}
\label{upc} The arguments are similar to those exposed for proving
the quasi convexity of the corridors. We give here only a sketch of
the proof. The horizontal deviation of an exceptional leaf with
respect to $\mathcal C$ depends linearly on the vertical variation
of the leaf (Lemma \ref{yoyoyo}). Thus, if a sufficiently large
segment of the leaf remains outside a sufficiently large horizontal
neighborhood of $\mathcal C$, the exponential separation of the
leaves implies that the horizontal distance between the leaf and
$\mathcal C$ exponentially increases with the vertical length of the
leaf. Assume now that the exceptional leaf considered is followed by
another exceptional one. The {\em strong} exponential separation
gives the same consequence: this second exceptional leaf does not go
back to $\mathcal C$ and the horizontal distance with respect to
$\mathcal C$ exponentially increases with its vertical length, as
soon as this length is sufficiently large enough. Here the arguments
are similar to those used for proving Lemmas \ref{stairs fibre} and
\ref{esperons}. Finally, if the exceptional leaf is followed by a
quasi geodesic in $\widehat{X}$, then the approximation by a stair
as was done before, yields the same conclusion. \hfill  $\Box$

\end{document}